*Type of Article* (Original Article)

# New Six Solutions to Solve Sixth Degree Polynomial Equation in General Forms by Relying on Radical Expressions

**Yassine Larbaoui[1]**

[1] *Department of Electrical Engineering, Université Hassan 1er, Settat, Morocco*
*y.larbaoui@uhp.ac.ma, yassine.larbaoui.uh1@gmail.com*

*Abstract -* This paper presents new six solutions for sixth degree polynomial equation in general forms basing on new theorems, where the possibility to calculate the six roots of any sixth degree equation nearly simultaneously. The proposed roots for sixth degree polynomials in this paper are structured basing on new proposed solutions for quartic polynomial equations, which we developed in order to reduce the expression of any sixth degree polynomial to an expression of fourth degree polynomial.

*Keywords -* new six solutions, new theorems, sixth degree polynomial, solving sixth degree equation.

## 1. Introduction

Polynomial equations are involved in different fields of science, computation and engineering, where the complexity of equations are dependent on the degrees of polynomials and their implicated coefficients. Therefore, solving polynomial equations with high degrees using formulary solutions is significantly complicated, where research efforts to find formulary roots are invested without fundamental methodologies to be extended and projected on polynomial equations with degrees higher than five.

Using radical expressions to solve polynomial equations with degrees higher than five becomes more complicated when considering general forms of polynomials, because of the necessity to eliminate radical expressions at the end of calculations and boot to solvable expressions in order to annihilate subparts of these polynomials.

Converging calculations on polynomial equations according to structured methodologies has less attention where efforts diverge to explore different paths of calculations. Therefore, shifting focus to the methodological aspect of analysing polynomials may have new insights.

Solving $n^{th}$ degree polynomial equations has been an enigmatic problem over hundreds of years, where many attempts concluded to the impossibility of elaborating universal solution formulas for polynomial equations with degrees equal or higher than fifth degree by using radicals. However, finding breakthroughs is a living hope for mathematicians who seek to solve quantic equations and above by using radical expressions.

Even though solving polynomial equations of fifth degree and higher by using radicals is center of focus, defining specific formulary solutions for fourth degree polynomial equations is still open to find different algebraic expressions, which may determine all four roots of any fourth degree polynomial in general form (nearly simultaneously) or bring different characteristics and new theorems.

This article presents the results of our research in mathematic science, where we present four formulary solutions for any fourth degree polynomial equation in general form with real coefficients, and we present six roots for any sixth degree polynomial equation in general form with real coefficients by using radical expressions. We attack the problem of solving polynomial equations of $n^{th}$ degree, where ($n = 4 \; and \; n = 6$), with a clear manner while using extendable logic, which enable us to elaborate different new theorems.

This paper is an introduction of a large research in mathematics to solve $n^{th}$ degree polynomial equations with new direct concepts and methods basing on projecting the presented work in this paper on other polynomial equations with degrees higher than four, which will be presented in other articles.

Lodovico Ferrari is attributed the credit of discovering a principal solution for fourth degree polynomial equation in 1540, but since his root expression required having a solution for cubic equation, which wasn't published yet by Gerolamo Cardano at that time, Lodovico Ferrari couldn't publish his discovery immediately. Ferrari's mentor, Gerolamo Cardano, in the book Ars Magna [1], published this discovered quartic solution along with the cubic solution.

The discovered cubic solution by Cardano [2] for third degree equations under the form $x^3 + px + q = 0$ helps Ferrari's solution to solve quartic equations by reducing them from the fourth degree to the second degree, but it doesn't directly help to properly define the four roots of any quartic equation. Therefore, when having a fourth degree polynomial equation in general form, it is necessary to conduct further calculations to determine its four roots while relying on Ferrari's solution. Cardano's method was also the base to solve particular forms of $n^{th}$ degree equations; such as by giving radical expressions under the form $\sqrt{na + \sqrt{b}} + \sqrt{na - \sqrt{b}}$ for $n = 2,3,4, \ldots$, etc. [3].

There are also other elaborated methods and proposed solutions for quartic equations such as Euler's solution [4], Galois's method [5], Descartes's method [6], Lagrange's method [7] and algebraic geometry [8].

Cubic solution was always essential as a base for further research to solve polynomial equations of fourth degree and above [9, 10], such as trying to prove that quantic equations do not accept quadratic expressions as solutions, which is discussed in [11-14].

A solution for third degree equation was first found in 1515 by Scipione del Ferro (1465–1526), for some specific cases defined by the values of coefficients. However, the official form of cubic solution, which is recognized as the base of further historical research on solving quartic equations and other specific forms of polynomial equations, is the published solution by Cardano.

There are also other recent researches dedicated to solve quartic equations where the used methods are based on expression reduction [15, 17], whereas other researches are relying on algorithms and numerical analysis to find the roots of polynomial equations with degrees higher than four [17, 18].

In the published paper by Tschirnhaus [19], he proposed an innovative method to solve polynomial equation $P_n(x)$ of $n^{th}$ degree by relying on its transformation into a reduced expression $Q_n(y)$ with fewer terms by extending the idea of Descartes; in which a polynomial of $n^{th}$ degree is reducible by removing its term in degree $(n - 1)$. The projection of this method on quantic equations is more detailed in [20].

There are some published papers treating resolvable quantic polynomials by relying on radical expressions [21], [22], where the description of specific characteristics which determine whether a quantic polynomial may accept roots with radical expressions or not. In addition, there are other published papers treating the resolvability of some polynomial equations with degrees higher than five by using factorization or by conditioning the forms of coefficients to be dependent on each other [23], [24], which do not solve neither quantic equations nor sixth degree equations in general forms.

Some polynomial equations of sixth degree in simple forms, such as $ax^6 + dx^3 + g = 0$, can be solved by factorizing into radicals, but other sixth degree equations in complete forms could not be solved over history [25], [26], until nowadays by proposing the presented solutions in this paper, which make the content of this article valuable.

The advantage of this paper is presenting, with details, the expressions of six roots for sixth degree polynomial equation in general form, where the possibility of calculating these roots nearly simultaneously by relying on our proposed solutions for quartic equations.

Because the content of this paper is original, and there are many new proposed formulas, mathematical expressions and theorems related in a scaling manner basing on extendable logic; every formula will be proved mathematically and used to build the rest of content, and we will go through them by logical analysis and deduction basing on structured development.

This paper is structured as follow: Section 2 where we present new four formulary solutions for fourth degree polynomial equation in general form enabling to calculate the roots of these equations nearly simultaneously. Section 3, where we present new six formulary solutions using radical expressions for sixth degree polynomial equation in general form where the coefficient of fifth degree part is different from zero. Section 4, where we present new six roots for sixth degree polynomial equations in general forms where the coefficient of fifth degree part is equal zero. Finally, Section 5 for conclusion.

## 2. Solutions and Theorems for Fourth Degree Polynomial Equation in General Form

This section presents new unified formulary solutions for fourth degree polynomial equation in general form.

### 2.1. *First Proposed Theorem*

In this subsection, we propose a new theorem to solve fourth degree polynomial equations that may be presented as shown in (eq.1). The expressions of proposed solutions are dependent on the value of the expression $\left(8\left(\frac{b}{a}\right)^3 - \frac{32cb}{a^2} + \frac{64d}{a}\right)$. We are proposing four solutions for $\left(8\left(\frac{b}{a}\right)^3 - \frac{32cb}{a^2} + \frac{64d}{a}\right) < 0$, four solutions for $\left(8\left(\frac{b}{a}\right)^3 - \frac{32cb}{a^2} + \frac{64d}{a}\right) > 0$ and four solutions for $\left(8\left(\frac{b}{a}\right)^3 - \frac{32cb}{a^2} + \frac{64d}{a}\right) = 0$.

The proposed solutions are expressed by using $y_{0,1}$, which is presented in (eq.17), and by using $P$, $Q$ and $R$ shown in (eq.6).

The proof of this theorem is presented in an independent subsection, because it is long and it contains the full expressions of proposed solutions along with details, in order to highlight the used logic to develop those solutions and to build the ground that we use to prove other proposed theorems in this paper.

### *Theorem 1*

A fourth degree polynomial equation under the expression (eq.1), where coefficients belong to the group of numbers $\mathbb{R}$, has four solutions:

$$ax^4 + bx^3 + cx^2 + dx + e = 0 \; where \; a \neq 0 \qquad (1)$$

If $\left(8\left(\frac{b}{a}\right)^3 - \frac{32cb}{a^2} + \frac{64d}{a}\right) < 0$, and by using the expressions of $y_{0,1}$ in (eq.17), P in (Eq.6) and Q in (Eq.6):

Solution 1: $S_{1,1}$ in expression (eq.35);

Solution 2: $S_{1,2}$ in expression (eq.36);

Solution 3: $S_{1,3}$ in expression (eq.37);

Solution 4: $S_{1,4}$ in expression (eq.38).

If $\left(8\left(\frac{b}{a}\right)^3 - \frac{32cb}{a^2} + \frac{64d}{a}\right) > 0$, and by using the expressions of $y_{0,1}$ in (eq.17), P in (Eq.6) and Q in (Eq.6):

Solution 1: $S_{2,1}$ in expression (eq.39);

Solution 2: $S_{2,2}$ in expression (eq.40);

Solution 3: $S_{2,3}$ in expression (eq.41);

Solution 4: $S_{2,4}$ in expression (eq.42).

If $\left(8\left(\frac{b}{a}\right)^3 - \frac{32cb}{a^2} + \frac{64d}{a}\right) = 0$, and by using the expressions of $y_{0,1}$ in (eq.28), P in (Eq.6) and Q in (Eq.6):

Solution 1: $S_{3,1}$ in expression (eq.43);

Solution 2: $S_{3,2}$ in expression (eq.44);

Solution 3: $S_{3,3}$ in expression (eq.45);

Solution 4: $S_{3,4}$ in expression (eq.46).

### 2.2. *Proof of Theorem 1*

By dividing the polynomial (eq.1) on the coefficient $a$, we have the next form:

$$x^4 + \frac{b}{a}x^3 + \frac{c}{a}x^2 + \frac{d}{a}x + \frac{e}{a} = 0 \; with \; a \neq 0 \quad (2)$$

We suppose that $x$ is expressed as shown in (eq.3):

$$x = \frac{\left(\frac{-b}{a}+y\right)}{4} \quad (3)$$

We replace $x$ with supposed expression in (eq.3) to reduce the form of presented polynomial in (eq.2). Thereby, we have the presented expression in (eq.4).

$$y^4 + y^2\left[-6\left(\frac{b}{a}\right)^2 + \frac{16c}{a}\right] + y\left[8\left(\frac{b}{a}\right)^3 - \frac{32cb}{a^2} + \frac{64d}{a}\right] - 3\left(\frac{b}{a}\right)^4 + \frac{16cb^2}{a^2} - \frac{64db}{a^2} + \frac{256e}{a} = 0 \quad (4)$$

To simplify the expression of shown polynomial equation in (eq.4), we replace the expressions of used coefficients as shown in (eq.5) where the values of those coefficients are defined in (eq.6).

$$y^4 + Py^2 + Qy + R = 0 \quad (5)$$

$$P = -6\left(\frac{b}{a}\right)^2 + \frac{16c}{a}; Q = 8\left(\frac{b}{a}\right)^3 - \frac{32cb}{a^2} + \frac{64d}{a}; R = -3\left(\frac{b}{a}\right)^4 + \frac{16cb^2}{a^2} - \frac{64db}{a^2} + \frac{256e}{a} \quad (6)$$

To solve the shown polynomial equation in (eq.5), we propose new expressions for the variable $y$; expressions (eq.7) and (eq.8):

$$\text{For } Q \leq 0: \quad y = \sqrt{y_0} + \sqrt{y_1} + \sqrt{y_2} \quad (7)$$

$$\text{For } Q \geq 0: \quad y = -\sqrt{y_0} - \sqrt{y_1} - \sqrt{y_2} \quad (8)$$

We propose the expressions (eq.9) and (eq.10) for $y_1$ and $y_2$ successively, in condition of $y_0 \neq 0$. Those expressions of $y_1$ and $y_2$ are based on quadratic solutions.

$$y_1 = -\frac{\frac{P}{2}+y_0}{2} + \sqrt{\left(\frac{\frac{P}{2}+y_0}{2}\right)^2 - \frac{Q^2}{64y_0}} \quad (9)$$

$$y_2 = -\frac{\frac{P}{2}+y_0}{2} - \sqrt{\left(\frac{\frac{P}{2}+y_0}{2}\right)^2 - \frac{Q^2}{64y_0}} \quad (10)$$

To reduce the expression of equation (eq.5) and find a way to solve it, we propose the following expressions for the coefficients $P$ and $Q$:

$$-2[y_0{}^2 + y_1{}^2 + y_2{}^2] = P \quad (11)$$

$$\text{For } Q \leq 0: \quad -8\sqrt{y_0}\sqrt{y_1}\sqrt{y_2} = Q \quad (12)$$

$$\text{For } Q \geq 0: \quad 8\sqrt{y_0}\sqrt{y_1}\sqrt{y_2} = Q \quad (13)$$

In the following calculation, we replace the variable $y$ with the expression (eq.7) where we suppose $Q < 0$, and we replace P and Q with their shown expressions in (eq.11) and (eq.12):

$$y^4 + Py^2 + Qy + R = -\left[\sqrt{y_0{}^4} + \sqrt{y_1{}^4} + \sqrt{y_2{}^4}\right] + 2\left[\sqrt{y_0{}^2}\sqrt{y_1{}^2} + \sqrt{y_0{}^2}\sqrt{y_2{}^2} + \sqrt{y_1{}^2}\sqrt{y_2{}^2}\right] + R$$

$$= -\left[\sqrt{y_0{}^2} + \sqrt{y_1{}^2} + \sqrt{y_2{}^2}\right]^2 + 4\left[\sqrt{y_0}{}^2\sqrt{y_1}{}^2 + \sqrt{y_0}{}^2\sqrt{y_2}{}^2 + \sqrt{y_1}{}^2\sqrt{y_2}{}^2\right] + R$$

$$= 4\left[-y_0\left(\frac{P}{2}+y_0\right)+\frac{Q^2}{64y_0}\right]+R-\frac{P^2}{4}=0$$

$$-\left[\sqrt{{y_0}^2}+\sqrt{{y_1}^2}+\sqrt{{y_2}^2}\right]^2+4\left[\sqrt{{y_0}^2}\sqrt{{y_1}^2}+\sqrt{{y_0}^2}\sqrt{{y_2}^2}+\sqrt{{y_1}^2}\sqrt{{y_2}^2}\right]+R=0 \Rightarrow y_0^3+\frac{P}{2}y_0^2+\frac{P^2-4R}{16}y_0-\frac{Q^2}{64}=0 \quad (14)$$

To solve the resulted expression in (eq.14), we use Cardano's solution for third degree polynomial equations.

For $w^3+cw+d=0$, Cardano's solution is as follow:

$$w=\sqrt[3]{\frac{-d}{2}+\sqrt{\left(\frac{d}{2}\right)^2+\left(\frac{c}{3}\right)^3}}+\sqrt[3]{\frac{-d}{2}-\sqrt{\left(\frac{d}{2}\right)^2+\left(\frac{c}{3}\right)^3}} \quad (15)$$

For $y^3+by^2+cy+d=0$, we use the form $y=\frac{-b+w}{3}$, and we suppose $\mathrm{D}=27d+2b^3-9cb$ and $\mathrm{C}=9c-3b^2$ to express the cubic solution as follow:

$$y=\frac{-b}{3}+\frac{1}{3}\sqrt[3]{-\frac{\mathrm{D}}{2}+\sqrt{\left(\frac{\mathrm{D}}{2}\right)^2+\left(\frac{\mathrm{C}}{3}\right)^3}}+\frac{1}{3}\sqrt[3]{-\frac{\mathrm{D}}{2}-\sqrt{\left(\frac{\mathrm{D}}{2}\right)^2+\left(\frac{\mathrm{C}}{3}\right)^3}} \quad (16)$$

By using the expression (eq.16), the solutions of third degree polynomial equation shown in (eq.14) are $y_{0,1}$ in (eq.17), $y_{0,2}$ in (eq.18) and $y_{0,3}$ in (eq.19), where $P^i=\frac{P}{2}$, $R^i=\frac{-27Q^2-2P^3+72PR}{64}$ and $Q^i=-\frac{3P^2+36R}{16}$

$$y_{0,1}=-\frac{P^i}{3}+\frac{1}{3}\sqrt[3]{-\frac{R^i}{2}+\sqrt{\left(\frac{R^i}{2}\right)^2+\left(\frac{Q^i}{3}\right)^3}}+\frac{1}{3}\sqrt[3]{-\frac{R^i}{2}-\sqrt{\left(\frac{R^i}{2}\right)^2+\left(\frac{Q^i}{3}\right)^3}} \quad (17)$$

In condition of $y_{0,1}\neq 0$, $y_{0,2}$ and $y_{0,3}$ are as follow:

$$y_{0,2}=-\frac{\frac{P}{2}+y_{0,1}}{2}+\sqrt{\left(\frac{\frac{P}{2}+y_{0,1}}{2}\right)^2-\frac{Q^2}{64y_{0,1}}} \quad (18)$$

$$y_{0,3}=-\frac{\frac{P}{2}+y_{0,1}}{2}-\sqrt{\left(\frac{\frac{P}{2}+y_{0,1}}{2}\right)^2-\frac{Q^2}{64y_{0,1}}} \quad (19)$$

We deduce that when $y_0$ takes the value $y_{0,1}$, the value of $y_{0,2}$ is equal to the shown value of $y_1$ in (eq.9), and the value of $y_{0,3}$ is equal to the shown value of $y_2$ in (eq.10).

There are three other possible expressions for $y$ which respect the proposition $-8\sqrt{y_0}\sqrt{y_1}\sqrt{y_2}=Q$ when $Q\leq 0$, and they give the same results of calculations toward having the shown third degree polynomial in (eq.14). Thereby, they give the same values for roots $y_{0,1},y_{0,2}$ and $y_{0,3}$. These three expressions are $y=-\sqrt{y_0}-\sqrt{y_1}+\sqrt{y_2}$, $y=-\sqrt{y_0}+\sqrt{y_1}-\sqrt{y_2}$ and $y=\sqrt{y_0}-\sqrt{y_1}-\sqrt{y_2}$.

By using the expressions in (eq.6) and (eq.17), the solutions of presented fourth degree polynomial equation in (eq.5) when $\left(8\left(\frac{b}{a}\right)^3-\frac{32cb}{a^2}+\frac{64d}{a}\right)<0$ are as shown in (eq.20), (eq.21), (eq.22) and (eq.23).

Solution 1: $\quad s_{1,1}=\sqrt{y_{0,1}}+\sqrt{-\frac{\frac{P}{2}+y_{0,1}}{2}+\sqrt{\left(\frac{\frac{P}{2}+y_{0,1}}{2}\right)^2-\frac{Q^2}{64y_{0,1}}}}+\sqrt{-\frac{\frac{P}{2}+y_{0,1}}{2}-\sqrt{\left(\frac{\frac{P}{2}+y_{0,1}}{2}\right)^2-\frac{Q^2}{64y_{0,1}}}} \quad (20)$

Solution 2: $\quad s_{1,2} = -\sqrt{y_{0,1}} - \sqrt{-\frac{\frac{P}{2}+y_{0,1}}{2} + \sqrt{\left(\frac{\frac{P}{2}+y_{0,1}}{2}\right)^2 - \frac{Q^2}{64y_{0,1}}}} + \sqrt{-\frac{\frac{P}{2}+y_{0,1}}{2} - \sqrt{\left(\frac{\frac{P}{2}+y_{0,1}}{2}\right)^2 - \frac{Q^2}{64y_{0,1}}}}$  (21)

Solution 3: $\quad s_{1,3} = -\sqrt{y_{0,1}} + \sqrt{-\frac{\frac{P}{2}+y_{0,1}}{2} + \sqrt{\left(\frac{\frac{P}{2}+y_{0,1}}{2}\right)^2 - \frac{Q^2}{64y_{0,1}}}} - \sqrt{-\frac{\frac{P}{2}+y_{0,1}}{2} - \sqrt{\left(\frac{\frac{P}{2}+y_{0,1}}{2}\right)^2 - \frac{Q^2}{64y_{0,1}}}}$  (22)

Solution 4: $\quad s_{1,4} = \sqrt{y_{0,1}} - \sqrt{-\frac{\frac{P}{2}+y_{0,1}}{2} + \sqrt{\left(\frac{\frac{P}{2}+y_{0,1}}{2}\right)^2 - \frac{Q^2}{64y_{0,1}}}} - \sqrt{-\frac{\frac{P}{2}+y_{0,1}}{2} - \sqrt{\left(\frac{\frac{P}{2}+y_{0,1}}{2}\right)^2 - \frac{Q^2}{64y_{0,1}}}}$  (23)

There are three other possible expressions for $y$ which respect the proposition $8\sqrt{y_0}\sqrt{y_1}\sqrt{y_2} = Q$ when $Q \geq 0$, and they give the same third degree polynomial shown in (eq.14) after calculations. These three expressions are $y = -\sqrt{y_0} + \sqrt{y_1} + \sqrt{y_2}$, $y = \sqrt{y_0} - \sqrt{y_1} + \sqrt{y_2}$ and $y = \sqrt{y_0} + \sqrt{y_1} - \sqrt{y_2}$.

By using the shown expressions in (eq.6) and (eq.17), the solutions of presented fourth degree polynomial equation in (eq.5) when $\left(8\left(\frac{b}{a}\right)^3 - \frac{32cb}{a^2} + \frac{64d}{a}\right) > 0$ are as shown in (eq.24), (eq.25), (eq.26) and (eq.27).

Solution 1: $\quad s_{2,1} = -\sqrt{y_{0,1}} - \sqrt{-\frac{\frac{P}{2}+y_{0,1}}{2} + \sqrt{\left(\frac{\frac{P}{2}+y_{0,1}}{2}\right)^2 - \frac{Q^2}{64y_{0,1}}}} - \sqrt{-\frac{\frac{P}{2}+y_{0,1}}{2} - \sqrt{\left(\frac{\frac{P}{2}+y_{0,1}}{2}\right)^2 - \frac{Q^2}{64y_{0,1}}}}$  (24)

Solution 2: $\quad s_{2,2} = -\sqrt{y_{0,1}} + \sqrt{-\frac{\frac{P}{2}+y_{0,1}}{2} + \sqrt{\left(\frac{\frac{P}{2}+y_{0,1}}{2}\right)^2 - \frac{Q^2}{64y_{0,1}}}} + \sqrt{-\frac{\frac{P}{2}+y_{0,1}}{2} - \sqrt{\left(\frac{\frac{P}{2}+y_{0,1}}{2}\right)^2 - \frac{Q^2}{64y_{0,1}}}}$  (25)

Solution 3: $\quad s_{2,3} = \sqrt{y_{0,1}} - \sqrt{-\frac{\frac{P}{2}+y_{0,1}}{2} + \sqrt{\left(\frac{\frac{P}{2}+y_{0,1}}{2}\right)^2 - \frac{Q^2}{64y_{0,1}}}} + \sqrt{-\frac{\frac{P}{2}+y_{0,1}}{2} - \sqrt{\left(\frac{\frac{P}{2}+y_{0,1}}{2}\right)^2 - \frac{Q^2}{64y_{0,1}}}}$  (26)

Solution 4: $\quad s_{2,4} = \sqrt{y_{0,1}} + \sqrt{-\frac{\frac{P}{2}+y_{0,1}}{2} + \sqrt{\left(\frac{\frac{P}{2}+y_{0,1}}{2}\right)^2 - \frac{Q^2}{64y_{0,1}}}} - \sqrt{-\frac{\frac{P}{2}+y_{0,1}}{2} - \sqrt{\left(\frac{\frac{P}{2}+y_{0,1}}{2}\right)^2 - \frac{Q^2}{64y_{0,1}}}}$  (27)

Concerning $\left(8\left(\frac{b}{a}\right)^3 - \frac{32cb}{a^2} + \frac{64d}{a}\right) = 0$:

The expression of $y_{0,1}$ is as shown in (eq.28) where $P^1 = \frac{P}{2}$, $R^1 = \frac{-2P^2+72PR}{64}$ and $Q^1 = -\frac{3P^2+36R}{16}$, whereas $y_{0,2}$ and $y_{0,3}$ are as shown in (eq.29) and (eq.30).

$$y_{0,1} = -\frac{P^1}{3} + \frac{1}{3}\sqrt[3]{-\frac{R^1}{2} + \sqrt{\left(\frac{R^1}{2}\right)^2 + \left(\frac{Q^1}{3}\right)^3}} + \frac{1}{3}\sqrt[3]{-\frac{R^1}{2} - \sqrt{\left(\frac{R^1}{2}\right)^2 + \left(\frac{Q^1}{3}\right)^3}} \quad (28)$$

$$y_{0,2} = -\frac{\frac{P}{2}+y_{0,1}}{2} + \sqrt{\left(\frac{\frac{P}{2}+y_{0,1}}{2}\right)^2 - \frac{Q^2}{64y_{0,1}}} = 0 \; or \; y_{0,2} = -\left(\frac{P}{2} + y_{0,2}\right) \quad (29)$$

$$y_{0,3} = -\frac{\frac{P}{2}+y_{0,1}}{2} - \sqrt{\left(\frac{\frac{P}{2}+y_{0,1}}{2}\right)^2 - \frac{Q^2}{64y_{0,1}}} = -\left(\frac{P}{2}+y_{0,2}\right) \ or \ y_{0,3} = 0 \quad (30)$$

Because of having the expressions $y_{0,2} = 0$ or $y_{0,3} = 0$, and having the intersection between the forms (eq.7) and (eq.8) for Q=0 $\left(y = \sqrt{y_0} + \sqrt{y_1} + \sqrt{y_2} \ for \ Q \leq 0 \ and \ y = -\sqrt{y_0} - \sqrt{y_1} - \sqrt{y_2} \ for \ Q \geq 0\right)$, there are four solutions for the polynomial equation shown in (eq.5) when Q=0 and they are as shown in (eq.31), (eq.32), (eq.33) and (eq.34).

Solution 1: $\quad s_{3,1} = \sqrt{y_{0,1}} + \sqrt{-\left(\frac{P}{2}+y_{0,1}\right)} \quad (31)$

Solution 2: $\quad s_{3,2} = -\sqrt{y_{0,1}} - \sqrt{-\left(\frac{P}{2}+y_{0,1}\right)} \quad (32)$

Solution 3: $\quad s_{3,3} = -\sqrt{y_{0,1}} + \sqrt{-\left(\frac{P}{2}+y_{0,1}\right)} \quad (33)$

Solution 4: $\quad s_{3,4} = \sqrt{y_{0,1}} - \sqrt{-\left(\frac{P}{2}+y_{0,1}\right)} \quad (34)$

When we give the value $y_{0,1}$ in (eq.17) to $y_0$, the values of $y_{0,2}$ in (eq.18) and $y_{0,3}$ in (eq.19) are equal to the shown values of $y_1$ in (eq.9) and $(y_3)$ in (eq.10) respectively. Thereby, even when we replace the value of $y_{0,1}$ in the expressions of proposed solutions by the values of $y_{0,2}$ or $y_{0,3}$, the results are only redundancies of proposed solutions, because the value of $P$ in the precedent expressions and in the proposed solutions is as follow:

$$\frac{P}{2} = -\left(\sqrt{y_0^2} + \sqrt{y_1^2} + \sqrt{y_2^2}\right) = -\left(\sqrt{y_{0,1}^2} + \sqrt{y_{0,2}^2} + \sqrt{y_{0,3}^2}\right)$$

In order to solve the polynomial equation shown in (eq.2), we use the expression $x = \frac{-\frac{b}{a}+y}{4}$ where $y$ is the unknown variable in polynomial equation (eq.5). By using expressions (eq.6) and (eq.17) for $\left(8\left(\frac{b}{a}\right)^3 - \frac{32cb}{a^2} + \frac{64d}{a}\right) < 0$ , the solutions for equation (eq.1) are as shown in (eq.35), (eq.36), (eq.37) and (eq.38).

Solution 1: $\quad S_{1,1} = -\frac{b}{4a} + \frac{1}{4}\sqrt{y_{0,1}} + \frac{1}{4}\sqrt{-\frac{\frac{P}{2}+y_{0,1}}{2} + \sqrt{\left(\frac{\frac{P}{2}+y_{0,1}}{2}\right)^2 - \frac{Q^2}{64y_{0,1}}}} + \frac{1}{4}\sqrt{-\frac{\frac{P}{2}+y_{0,1}}{2} - \sqrt{\left(\frac{\frac{P}{2}+y_{0,1}}{2}\right)^2 - \frac{Q^2}{64y_{0,1}}}} \quad (35)$

Solution 2: $\quad S_{1,2} = -\frac{b}{4a} - \frac{1}{4}\sqrt{y_{0,1}} - \frac{1}{4}\sqrt{-\frac{\frac{P}{2}+y_{0,1}}{2} + \sqrt{\left(\frac{\frac{P}{2}+y_{0,1}}{2}\right)^2 - \frac{Q^2}{64y_{0,1}}}} + \frac{1}{4}\sqrt{-\frac{\frac{P}{2}+y_{0,1}}{2} - \sqrt{\left(\frac{\frac{P}{2}+y_{0,1}}{2}\right)^2 - \frac{Q^2}{64y_{0,1}}}} \quad (36)$

Solution 3: $\quad S_{1,3} = -\frac{b}{4a} - \frac{1}{4}\sqrt{y_{0,1}} + \frac{1}{4}\sqrt{-\frac{\frac{P}{2}+y_{0,1}}{2} + \sqrt{\left(\frac{\frac{P}{2}+y_{0,1}}{2}\right)^2 - \frac{Q^2}{64y_{0,1}}}} - \frac{1}{4}\sqrt{-\frac{\frac{P}{2}+y_{0,1}}{2} - \sqrt{\left(\frac{\frac{P}{2}+y_{0,1}}{2}\right)^2 - \frac{Q^2}{64y_{0,1}}}} \quad (37)$

Solution 4: $\quad S_{1,4} = -\frac{b}{4a} + \frac{1}{4}\sqrt{y_{0,1}} - \frac{1}{4}\sqrt{-\frac{\frac{P}{2}+y_{0,1}}{2} + \sqrt{\left(\frac{\frac{P}{2}+y_{0,1}}{2}\right)^2 - \frac{Q^2}{64y_{0,1}}}} - \frac{1}{4}\sqrt{-\frac{\frac{P}{2}+y_{0,1}}{2} - \sqrt{\left(\frac{\frac{P}{2}+y_{0,1}}{2}\right)^2 - \frac{Q^2}{64y_{0,1}}}} \quad (38)$

By using the expression $x = \frac{-\frac{b}{a}+y}{4}$ while relying on expressions (eq.6) and (eq.17) for $\left(8\left(\frac{b}{a}\right)^3 - \frac{32cb}{a^2} + \frac{64d}{a}\right) > 0$, the proposed solutions for equation (eq.1) are as shown in (eq.39), (eq.40), (eq.41) and (eq.42).

Solution 1: $S_{2,1} = -\frac{b}{4a} - \frac{1}{4}\sqrt{y_{0,1}} - \frac{1}{4}\sqrt{-\frac{\frac{P}{2}+y_{0,1}}{2} + \sqrt{\left(\frac{\frac{P}{2}+y_{0,1}}{2}\right)^2 - \frac{Q^2}{64y_{0,1}}}} - \frac{1}{4}\sqrt{-\frac{\frac{P}{2}+y_{0,1}}{2} - \sqrt{\left(\frac{\frac{P}{2}+y_{0,1}}{2}\right)^2 - \frac{Q^2}{64y_{0,1}}}}$ (39)

Solution 2: $S_{2,2} = -\frac{b}{4a} - \frac{1}{4}\sqrt{y_{0,1}} + \frac{1}{4}\sqrt{-\frac{\frac{P}{2}+y_{0,1}}{2} + \sqrt{\left(\frac{\frac{P}{2}+y_{0,1}}{2}\right)^2 - \frac{Q^2}{64y_{0,1}}}} + \frac{1}{4}\sqrt{-\frac{\frac{P}{2}+y_{0,1}}{2} - \sqrt{\left(\frac{\frac{P}{2}+y_{0,1}}{2}\right)^2 - \frac{Q^2}{64y_{0,1}}}}$ (40)

Solution 3: $S_{2,3} = -\frac{b}{4a} + \frac{1}{4}\sqrt{y_{0,1}} - \frac{1}{4}\sqrt{-\frac{\frac{P}{2}+y_{0,1}}{2} + \sqrt{\left(\frac{\frac{P}{2}+y_{0,1}}{2}\right)^2 - \frac{Q^2}{64y_{0,1}}}} + \frac{1}{4}\sqrt{-\frac{\frac{P}{2}+y_{0,1}}{2} - \sqrt{\left(\frac{\frac{P}{2}+y_{0,1}}{2}\right)^2 - \frac{Q^2}{64y_{0,1}}}}$ (41)

Solution 4: $S_{2,4} = -\frac{b}{4a} + \frac{1}{4}\sqrt{y_{0,1}} + \frac{1}{4}\sqrt{-\frac{\frac{P}{2}+y_{0,1}}{2} + \sqrt{\left(\frac{\frac{P}{2}+y_{0,1}}{2}\right)^2 - \frac{Q^2}{64y_{0,1}}}} - \frac{1}{4}\sqrt{-\frac{\frac{P}{2}+y_{0,1}}{2} - \sqrt{\left(\frac{\frac{P}{2}+y_{0,1}}{2}\right)^2 - \frac{Q^2}{64y_{0,1}}}}$ (42)

By using the expression $x = -\frac{\frac{b}{a}+y}{4}$ while relying on expressions (eq.6) and (eq.28) for $\left(8\left(\frac{b}{a}\right)^3 - \frac{32cb}{a^2} + \frac{64d}{a}\right) = 0$, the proposed solutions for equation (eq.1) are as shown in (eq.43), (eq.44), (eq.45) and (eq.46).

Solution 1: $S_{3,1} = -\frac{b}{4a} + \frac{1}{4}\sqrt{y_{0,1}} + \frac{1}{4}\sqrt{-\left(\frac{P}{2}+y_{0,1}\right)}$ (43)

Solution 2: $S_{3,2} = -\frac{b}{4a} - \frac{1}{4}\sqrt{y_{0,1}} - \frac{1}{4}\sqrt{-\left(\frac{P}{2}+y_{0,1}\right)}$ (44)

Solution 3: $S_{3,3} = -\frac{b}{4a} - \frac{1}{4}\sqrt{y_{0,1}} + \frac{1}{4}\sqrt{-\left(\frac{P}{2}+y_{0,1}\right)}$ (45)

Solution 4: $S_{3,4} = -\frac{b}{4a} + \frac{1}{4}\sqrt{y_{0,1}} - \frac{1}{4}\sqrt{-\left(\frac{P}{2}+y_{0,1}\right)}$ (46)

## 3. New Six Solutions for Sixth Degree Polynomial Equation in General Form

In this section, we propose six new solutions for sixth degree polynomial equation in general form shown in (eq.47), where we rely on our proposed solutions for quartic polynomial equations in previous section to structure the expressions of proposed solutions for sixth degree polynomial equations. We extend the used logic in precedent theorem (Theorem 1) by projection on sixth degree equations to prove the expressions of developed roots.

### 3.1. Second Proposed Theorem

In this subsection, we present our second proposed theorem to introduce new six formulary solutions for sixth degree polynomial equation in general form shown in (eq.47), where coefficients belong to the group $\mathbb{R}$ whereas the coefficient of fifth degree part is different from zero. First, we divide the polynomial (eq.47) on coefficient A to reduce its expression to the simplified form shown in (eq.48) where coefficients are expressed as shown in (eq.49).

$$Ax^6 + Bx^5 + Cx^4 + Dx^3 + Ex^2 + Fx + G = 0 \text{ with } A \neq 0 \text{ and } B \neq 0 \quad (47)$$

$$x^6 + bx^5 + cx^4 + dx^3 + ex^2 + fx + g = 0 \text{ with } b \neq 0 \quad (48)$$

$$b = \frac{B}{A}; \ c = \frac{C}{A}; \ D = \frac{D}{A}; e = \frac{E}{A}; \ f = \frac{F}{A}; \ g = \frac{G}{A}; \quad (49)$$

$$z^4 + \Gamma_3 z^3 + \Gamma_2 z^2 + \Gamma_1 z + \Gamma_0 = 0 \quad (50)$$

**Theorem 2**

After reducing the form of sixth degree polynomial shown in (eq.47) to the presented form in (eq.48) where coefficients are as expressed in (eq.49); the sixth degree polynomial equation shown in (eq.48), where coefficients belong to the group of numbers $\mathbb{R}$, can be reduced to a fourth degree polynomial equation, which may be expressed as shown in (eq.50). The reduction from sixth degree polynomial to quartic polynomial is conducted by supposing $x = x_0 x_1 + x_0 x_2 + x_0 x_3 + x_1 x_2 + x_1 x_3 + x_2 x_3$, whereas supposing $z = (x_0 + x_1 + x_2 + x_3)$ is the solution for fourth degree polynomial equation in (eq.50) by using Theorem 1 and relying on the expression $x_3 = -\frac{\Gamma_3}{4}$. The variable $\Gamma_3$ is defined as shown in (eq.51) where $\alpha_3$ is presented in (eq.52) and $\Gamma_4$ is the solution for the polynomial equation (eq.53), which relies on the coefficients (eq.54), (eq.55), (eq.56) and (eq.57). The shown coefficients in (eq.54), (eq.55), (eq.56) and (eq.57) are expressed by using the constant $V$ which is presented in (eq.58). The coefficients $\Gamma_3, \Gamma_2, \Gamma_1$ and $\Gamma_0$ of quartic equation (eq.50), which is used to calculate $z$, are determined by using the shown expressions in (51), (eq.59), (eq.60) and (eq.61) while using calculated values of $\Gamma_4$ and $V$. As a result, we have twelve calculated values as potential solutions for sixth degree polynomial equation shown in (eq.48), where many of them are only redundancies of others, because there are only six official solutions to determine.

The twelve solutions to calculate for sixth degree equation (eq.48) are as shown in the groups (eq.99), (eq.100) and (eq.101). The proposed six values as official solutions for sixth degree polynomial equation shown in (eq.48) are as presented in (eq.102), (eq.103), (eq.104), (eq.105), (eq.106) and (eq.107).

$$\Gamma_3 = \frac{4\alpha_3}{b} + \Gamma_4 \quad (51)$$

$$\alpha_3 = -\frac{\frac{4\Gamma_4\left(f - \frac{d^2}{4b}\right)}{b}}{\frac{32f}{b^2} + \frac{40d^2}{b^3} - \frac{64cd}{b^2} + \frac{64e}{b}} \quad (52)$$

$$\lambda_3\Gamma_4^6 + \lambda_2\Gamma_4^4 + \lambda_1\Gamma_4^2 + \lambda_0 = 0 \quad (53)$$

$$\lambda_3 = -\frac{40960}{V^4 b^4} + \frac{16384}{V^3 b^3} - \frac{1536}{V^2 b^2} \quad (54)$$

$$\lambda_2 = -\frac{24576d}{V^2 b^4} + \frac{16384c}{V^2 b^3} + \frac{3072d}{Vb^3} - \frac{2048c}{Vb^2} + \frac{1024}{V} \quad (55)$$

$$\lambda_1 = -\frac{512d}{b} + \frac{1536f}{b^3} + \frac{28V^2f}{b} - \frac{7V^2d^2}{b^2} + \frac{96Vf}{b^2} - \frac{168d^2V}{b^3} + \frac{192cdV}{b^2} - \frac{192Ve}{b} - \frac{3456d^2}{b^4} + \frac{4096cd}{b^3} - \frac{1024e}{b^2} - \frac{1024c^2}{b^2} \quad (56)$$

$$\lambda_0 = -\frac{64V^2d^2}{b^4} + \frac{64cd^2V^2}{b^3} - \frac{64eV^2d}{b^2} + \frac{128V^2g}{b} + \frac{192V^2df}{b^3} - \frac{128V^2cf}{b^2} \quad (57)$$

$$V = -\frac{\frac{32f}{b^2} + \frac{40d^2}{b^3} - \frac{64cd}{b^2} + \frac{64e}{b}}{4\frac{\left(f - \frac{d^2}{4b}\right)}{b}} \quad (58)$$

$$\Gamma_2 = \frac{8\Gamma_4^2}{Vb} - \frac{6d}{b^2} + \frac{4c}{b} + \frac{\left(f - \frac{d^2}{4b}\right)V^2}{2b\Gamma_4^2} - \frac{8\Gamma_4^2}{V^2 b^2} \quad (59)$$

$$\Gamma_1 = \frac{5\Gamma_4^3}{Vb} + \frac{3Vd^2}{4b^2\Gamma_4} - \frac{6d\Gamma_4}{b^2} + \frac{4c\Gamma_4}{b} - \frac{dcV}{b^2\Gamma_4} + \frac{eV}{b\Gamma_4} - \frac{\Gamma_4^3}{4} - \frac{8\Gamma_4^3}{V^2 b^2} + \frac{f - \frac{d^2}{4b}}{4\Gamma_4 b}V^2 \quad (60)$$

$$\Gamma_0 = \frac{\Gamma_4^4}{2Vb} - \frac{V^2d^3}{16b^4\Gamma_4^2} + \frac{3Vd^2}{8b^3} - \frac{3d\Gamma_4^2}{4b^2} + \frac{c\Gamma_4^2}{2b} + \frac{cd^2V^2}{8b^3\Gamma_4^2} - \frac{cdV}{2b^2} + \frac{eV}{2b} - \frac{eV^2d}{4b^2\Gamma_4^2} + \frac{gV^2}{2b\Gamma_4^2} - \left(\frac{\Gamma_4^2}{4} + V^2\frac{f - \frac{d^2}{4b}}{4b\Gamma_4^2}\right)\left(\frac{\Gamma_4^2}{4} + \frac{8\Gamma_4^2}{V^2 b^2} - \frac{2\Gamma_4^2}{VB} + \frac{3d}{b^2} - \frac{2c}{b} - \frac{\left(f - \frac{d^2}{4b}\right)}{4b\Gamma_4^2}V^2\right)$$

$$(61)$$

### 3.2. Proof of Theorem 2

Considering the sixth degree polynomial equation shown in (eq.48), we propose the expression (eq.62) in order to reduce the form from sixth degree to a fourth degree polynomial. We also propose the expression (eq.63), which presents the solution form for quartic equation by extending the used logic and presented solutions in Theorem 1.

$$x = x_0 x_1 + x_0 x_2 + x_0 x_3 + x_1 x_2 + x_1 x_3 + x_2 x_3 \quad (62)$$

$$z = x_0 + x_1 + x_2 + x_3 \quad (63)$$

We replace $x$ with its proposed value in (eq.62), in order to end by calculations to the reduced form shown in (eq.50).

In the shown expressions in (eq.64), we rely on the use of $x_i$, $x_j$ and $x_k$ where $\{x_i, x_j, x_k\} \in \{x_0, x_1, x_2, x_3\}$ and $i \neq j \neq k$.

$$\alpha_1 = \sum_{i=0}^{i=3} x_i^2 \; ; \; \alpha_2 = \sum_{i \neq j} x_i^2 x_j^2 ; \; \alpha_3 = \sum_{i \neq j \neq k} x_i x_j x_k ; \; \alpha_4 = x_0 x_1 x_2 x_3 \quad (64)$$

$$x = \frac{[z^2 - \alpha_1]}{2} \quad (65)$$

$$x^2 = \alpha_2 + 2\alpha_3 z + 6\alpha_4 \quad (66)$$

$$x^3 = z^3 \alpha_3 + \frac{1}{2} z^2 [\alpha_2 + 6\alpha_4] - z\alpha_1 \alpha_3 - \frac{1}{2} [\alpha_2 + 6\alpha_4]\alpha_1 \quad (67)$$

$$x^4 = 4z^2 \alpha_3^2 + 4\alpha_3 z[\alpha_2 + 6\alpha_4] + [\alpha_2 + 6\alpha_4]^2 \quad (68)$$

$$x^5 = 2z^4 \alpha_3^2 + 2z^3 [\alpha_2 + 6\alpha_4]\alpha_3 + \frac{1}{2} z^2 [(\alpha_2 + 6\alpha_4)^2 - 4\alpha_3^2 \alpha_1] - 2z\alpha_3 \alpha_1 [\alpha_2 + 6\alpha_4] - \frac{1}{2} \alpha_1 [\alpha_2 + 6\alpha_4]^2$$

$$(69)$$

$$x^6 = 8z^3 \alpha_3^3 + 12z^2 [\alpha_2 + 6\alpha_4]\alpha_3^2 + 6z[(\alpha_2 + 6\alpha_4)]^2 \alpha_3 + [\alpha_2 + 6\alpha_4]^3 \quad (70)$$

$$\gamma_4 z^4 + \gamma_3 z^3 + \gamma_2 z^2 + \gamma_1 z + \gamma_0 = 0 \quad (71)$$

We use the expressions of $\{\alpha_1, \alpha_2, \alpha_3, \alpha_4\}$ in (eq.64), $x$ in (eq.65), $x^2$ in (eq.66), $x^3$ in (eq.67), $x^4$ in (eq.68), $x^5$ in (eq.69) and $x^6$ in (eq.70), to have the fourth degree polynomial shown in (eq.71) where the values of coefficients are as follow:

$$\gamma_4 = 2b\alpha_3^2$$

$$\gamma_3 = 8\alpha_3^3 + d\alpha_3 + 2b[\alpha_2 + 6\alpha_4]\alpha_3$$

$$\gamma_2 = 12[\alpha_2 + 6\alpha_4]\alpha_3^2 + \frac{1}{2} b[(\alpha_2 + 6\alpha_4)^2 - 4\alpha_3^2 \alpha_1] + 4c\alpha_3^2 + \frac{1}{2} d[\alpha_2 + 6\alpha_4] + \frac{1}{2} f$$

$$\gamma_1 = 6[(\alpha_2 + 6\alpha_4)]^2 \alpha_3 - 2b\alpha_3 \alpha_1 [\alpha_2 + 6\alpha_4] + 4c[\alpha_2 + 6\alpha_4]\alpha_3 - d\alpha_1 \alpha_3 + 2e\alpha_3$$

$$\gamma_0 = [\alpha_2 + 6\alpha_4]^3 - \frac{1}{2} b\alpha_1 [\alpha_2 + 6\alpha_4]^2 + c[\alpha_2 + 6\alpha_4]^2 - \frac{1}{2} d[\alpha_2 + 6\alpha_4]\alpha_1 + e[\alpha_2 + 6\alpha_4] - \frac{1}{2} f\alpha_1 + g$$

We divide the polynomial equation shown in (eq.71) on $\gamma_4$ to simplify its expression. As a result, we have the shown equation in (eq.50) where the values of coefficients are as follow:

$$\Gamma_4 = \frac{[\alpha_2 + 6\alpha_4] + \frac{d}{2b}}{\alpha_3} \Rightarrow \alpha_2 = \alpha_3 \Gamma_4 - \frac{d}{2b} - 6\alpha_4$$

$$\Gamma_3 = \frac{4\alpha_3}{b} + \Gamma_4$$

$$\Gamma_2 = \frac{6\left[\Gamma_4\alpha_3 - \frac{d}{2b}\right]}{b} + \frac{\Gamma_4^2}{4} - \alpha_1 + \frac{2c}{b} + \frac{f - \frac{d^2}{4b}}{4b\alpha_3^2}$$

$$\Gamma_1 = \frac{3\left[\Gamma_4\alpha_3 - \frac{d}{2b}\right]^2}{b\alpha_3} - \Gamma_4\alpha_1 + \frac{2c\left[\alpha_3\Gamma_4 - \frac{d}{2b}\right]}{b\alpha_3} + \frac{e}{b\alpha_3}$$

$$\Gamma_0 = \frac{\left[\Gamma_4\alpha_3 - \frac{d}{2b}\right]^3}{2b\alpha_3^2} - \frac{1}{4}\alpha_1\Gamma_4^2 + \frac{c\left[\alpha_3\Gamma_4 - \frac{d}{2b}\right]^2}{2b\alpha_3^2} + \frac{e\left[\alpha_3\Gamma_4 - \frac{d}{2b}\right]}{2b\alpha_3^2} + \frac{g - \frac{f\alpha_1}{2} + \frac{d^2\alpha_1}{8b}}{2b\alpha_3^2}$$

From precedent section (section 2), we ended with solutions expressed as $y = y_0 + y_1 + y_2$ for fourth degree polynomial equation in simple form shown in (eq.5), whereas the solution for fourth degree polynomial equation in complete form is expressed as $z = -\frac{b}{4a} + \frac{1}{4}y_0 + \frac{1}{4}y_1 + \frac{1}{4}y_2$.

We replace $z$ with $\frac{-\Gamma_3 + y}{4}$ in order to reduce the form of quartic equation from expression (eq.50) to expression (eq.72), where the values of coefficients are as shown in (eq.73).

$$y^4 + Py^2 + Qy + R = 0 \quad (72)$$

$$P = -6\Gamma_3^2 + 16\Gamma_2; Q = 8\Gamma_3^3 - 32\Gamma_2\Gamma_3 + 64\Gamma_1; R = -3\Gamma_3^4 + 16\Gamma_2\Gamma_3^2 - 64\Gamma_1\Gamma_3 + 256\Gamma_0 \quad (73)$$

Concerning the fourth degree polynomial equation in (eq.50), where $z$ is as shown in (eq.63), the principal proposed expressions for the solutions are $z = -\frac{\Gamma_3}{4} + \frac{1}{4}\sqrt{y_0} + \frac{1}{4}\sqrt{y_1} + \frac{1}{4}\sqrt{y_2}$ when $Q \leq 0$ and $z = -\frac{\Gamma_3}{4} - \frac{1}{4}\sqrt{y_0} - \frac{1}{4}\sqrt{y_1} - \frac{1}{4}\sqrt{y_2}$ *when* $Q \geq 0$; where $x_3 = -\frac{\Gamma_3}{4}, x_0 = \pm\frac{1}{4}\sqrt{y_0}, x_1 = \pm\frac{1}{4}\sqrt{y_1}$ and $x_2 = \pm\frac{1}{4}\sqrt{y_2}$. These two principal expressions are sufficient to conduct the calculations of proof, and then generalize the results by using the other expressed forms of solutions in Theorem 1.

We replace $\Gamma_3, \Gamma_2, \Gamma_1$ *and* $\Gamma_0$ with their values in function of $\{\Gamma_4 , \alpha_1 , \alpha_3\}$ in order to have the expressions of $P$ in (eq.74), $Q$ in (eq.75) and $R$ in (eq.76).

$$P = -2\Gamma_4^2 - \frac{96\alpha_3^2}{b^2} + \frac{48\Gamma_4\alpha_3}{b} - \frac{48d}{b^2} - 16\alpha_1 + \frac{32c}{b} + 4\frac{\left(f - \frac{d^2}{4b}\right)}{b\alpha_3^2} \quad (74)$$

$$Q = \frac{512\alpha_3^3}{b^3} - \frac{384\alpha_3^2\Gamma_4}{b^2} + \frac{64\alpha_3\Gamma_4^2}{b} + \frac{384d\alpha_3}{b^3} + \frac{128\alpha_3\alpha_1}{b} - \frac{256c\alpha_3}{b^2} - \frac{32f}{b^2\alpha_3} - \frac{96\Gamma_4}{b^2} - 32\Gamma_4\alpha_1 + \frac{64c\Gamma_4}{b} - \frac{8\Gamma_4\left(f - \frac{d^2}{4b}\right)}{b\alpha_3^2} + \frac{56d^2}{b^3\alpha_3} - \frac{64cd}{b^2\alpha_3} + \frac{64e}{b\alpha_3}$$

$$(75)$$

$$R = -\frac{768\alpha_3^4}{b^4} + \Gamma_4^4 + \frac{768\Gamma_4\alpha_3^3}{b^3} + \frac{16\Gamma_4^3\alpha_3}{b} - \frac{224\Gamma_4^2\alpha_3^2}{b^2} - \frac{768d\alpha_3^2}{b^4} - \frac{256\alpha_3^2\alpha_1}{b^2} + \frac{512c\alpha_3^2}{b^3} + \frac{64f}{b^3} - \frac{48d\Gamma_4^2}{b^2} - 16\alpha_1\Gamma_4^2 + \frac{32c\Gamma_4^2}{b} + 4\frac{\left(f - \frac{d^2}{4b}\right)}{b\alpha_3^2}\Gamma_4^2 + \frac{384d\Gamma_4\alpha_3}{b^3} + \frac{128\Gamma_4\alpha_3\alpha_1}{b} - \frac{256c\Gamma_4\alpha_3}{b^2} + \frac{32f\Gamma_4}{b^2\alpha_3} + \frac{40d^2\Gamma_4}{b^3\alpha_3} - \frac{208d^2}{b^4} + \frac{256cd}{b^3} - \frac{256e}{b^2} - \frac{64cd\Gamma_4}{b^2\alpha_3} + \frac{64e\Gamma_4}{b\alpha_3} - \frac{16d^3}{b^4\alpha_3^2} + \frac{32cd^2}{b^3\alpha_3^2} - \frac{64ed}{b^2\alpha_3^2} + 128\frac{\left(g - \frac{f\alpha_1}{2} + \frac{d^2\alpha_1}{8b}\right)}{b\alpha_3^2}$$

$$(76)$$

By using the expressions (eq.11), (eq.12), (eq.13) and (eq.14) from the proof of first proposed theorem (Theorem 1), the values of $P$, $Q$ and $R$ are as shown in (eq.77), (eq.78) and (eq.79) successively.

$$P = -2[(4x_0)^2 + (4x_1)^2 + (4x_2)^2] \;\Rightarrow\; P = -32\left[\alpha_1 - \frac{\Gamma_3^2}{16}\right] = 2\Gamma_4^2 + \frac{16\Gamma_4\alpha_3}{b} + \frac{32\alpha_3^2}{b^2} - 32\alpha_1 \quad (77)$$

$$Q = -8(4x_0)(4x_1)(4x_2) \;\Rightarrow\; Q = -\frac{8(4x_0)(4x_1)(4x_2)(4x_3)}{4x_3} = \frac{2048\alpha_4}{\Gamma_3} = \frac{2048\alpha_4}{\frac{4\alpha_3}{b}+\Gamma_4} \quad (78)$$

$$R = [(4x_0)^2 + (4x_1)^2 + (4x_2)^2]^2 - 4[(4x_0)^2(4x_1)^2 + (4x_0)^2(4x_2)^2 + (4x_1)^2(4x_2)^2]$$

$$\Rightarrow R = 256\left[\alpha_1 - \frac{\Gamma_3^2}{16}\right]^2 - 1024\left[\alpha_2 - \frac{\Gamma_3^2}{16}\left(\alpha_1 - \frac{\Gamma_3^2}{16}\right)\right]$$

$$\Rightarrow R = -\frac{768\alpha_3^4}{b^4} - 3\Gamma_4^4 - \frac{768\alpha_3^3\Gamma_4}{b^3} - \frac{48\Gamma_4^3\alpha_3}{b} - \frac{288\Gamma_4^2\alpha_3^2}{b^2} + \frac{512\alpha_3^2\alpha_1}{b^2} + 32\alpha_1\Gamma_4^2 + \frac{256\alpha_1\alpha_3\Gamma_4}{b} + 256\alpha_1^2 - 1024\alpha_2$$

$$(79)$$

We have a group of four variables $\{\alpha_1, \alpha_2, \alpha_3, \alpha_4\}$, whereas we have a group of only three equations to solve $\{P, Q, R\}$ where all of them are dependent on the value of $\Gamma_4$. Thereby, the next step is about using the appropriate logic of analysis and calculation to find the value of $\Gamma_4$ while taking advantage of the fact that having a group of four variables enables us to solve four equations.

In order to reduce the expression of $R$ in (eq.76), and find a way to determine the value of $\Gamma_4$, we suppose that $\left(\frac{4\Gamma_4^2\left(f-\frac{d^2}{4b}\right)}{b\alpha_3^2} + \frac{32f\Gamma_4}{b^2\alpha_3} + \frac{40d^2\Gamma_4}{b^3\alpha_3} - \frac{64cd\Gamma_4}{b^2\alpha_3} + \frac{64e\Gamma_4}{b\alpha_3}\right) = 0$ where $\frac{\Gamma_4}{\alpha_3} \neq 0$. As a result, we have the shown expression in (eq.80).

$$\frac{\Gamma_4}{\alpha_3} = V = -\frac{\frac{32f}{b^2} + \frac{40d^2}{b^3} - \frac{64cd}{b^2} + \frac{64e}{b}}{\frac{4\left(f-\frac{d^2}{4b}\right)}{b}} \quad (80)$$

From expression (eq.80), we can see that $\frac{\Gamma_4}{\alpha_3}$ has a constant value. Therefore, in the rest of calculation, we will replace $\frac{\Gamma_4}{\alpha_3}$ by the constant $V$, which is shown in (eq.58).

From precedent calculations of dividing the polynomial equation shown in (eq.71) on $\gamma_4$ to simplify its expression, along using the expression of $\Gamma_4$, we have $\alpha_2 = \left(\alpha_3\Gamma_4 - \frac{d}{2b} - 6\alpha_4\right)$.

We have the resulted equation in (eq.81) by using the expressions of $P$ in (eq.74) and (eq.77), which we use to define the shown value of $\alpha_1$ in (eq.82).

$$-4\Gamma_4^2 - \frac{128\Gamma_4^2}{V^2b^2} + \frac{32\Gamma_4^2}{Vb} - \frac{48d}{b^2} + 16\alpha_1 + \frac{32c}{b} + \frac{4V^2\left(f-\frac{d^2}{4b}\right)}{\Gamma_4^2b} = 0 \quad (81)$$

$$\alpha_1 = \frac{\Gamma_4^4 + \frac{32\Gamma_4^4}{V^2b^2} - \frac{8\Gamma_4^4}{Vb} + \frac{12d\Gamma_4^2}{b^2} - \frac{8c\Gamma_4^2}{b} - \frac{V^2\left(f-\frac{d^2}{4b}\right)}{b}}{4\Gamma_4^2} \quad (82)$$

We have the presented polynomial equation in (eq.83) by using the expressions of $R$ in (eq.76) and (eq.79).

$$4\Gamma_4^4 + \frac{1536\Gamma_4\alpha_3^3}{b^3} + \frac{64\Gamma_4^3\alpha_3}{b} + \frac{64\Gamma_4^2\alpha_3^2}{b^2} - \frac{768d\alpha_3^2}{b^4} - \frac{768\alpha_3^2\alpha_1}{b^2} + \frac{512c\alpha_3^3}{b^3} + \frac{64f}{b^3} - \frac{48d\Gamma_4^2}{b^2}$$
$$-48\alpha_1\Gamma_4^2 + \frac{32c\Gamma_4^2}{b} + \frac{384d\Gamma_4\alpha_3}{b^3} - \frac{128\Gamma_4\alpha_3\alpha_1}{b} - \frac{256c\Gamma_4\alpha_3}{b^2} - \frac{208d^2}{b^4} + \frac{256cd}{b^3} - \frac{256e}{b^2} - \frac{16d^3}{b^4\alpha_3} + \frac{32cd^2}{b^3\alpha_3^2} -$$
$$\frac{64ed}{b^2\alpha_3^2} - 256\alpha_1^2 + 1024\alpha_2 + 128\frac{\left(g - \frac{f\alpha_1}{2} + \frac{d^2\alpha_1}{8b}\right)}{b\alpha_3^2} = 0$$

(83)

We replace $\frac{\Gamma_4}{\alpha_3}$ with its shown expression in (eq.80) to pass from equation (eq.83) to equation (eq.84).

$$4\Gamma_4^4 + \frac{1536\Gamma_4^4}{V^3b^3} + \frac{64\Gamma_4^4}{Vb} + \frac{64\Gamma_4^4}{V^2b^2} - \frac{768d\Gamma_4^2}{V^2b^4} - \frac{768\Gamma_4^2\alpha_1}{V^2b^2} + \frac{512c\Gamma_4^2}{V^2b^3} + \frac{64f}{b^3} - \frac{48d\Gamma_4^2}{b^2} - 48\alpha_1\Gamma_4^2 + \frac{32c\Gamma_4^2}{b} +$$
$$\frac{384d\Gamma_4^2}{Vb^3} - \frac{128\Gamma_4^2\alpha_1}{Vb} - \frac{256c\Gamma_4^2}{Vb^2} - \frac{208d^2}{b^4} + \frac{256cd}{b^3} - \frac{256e}{b^2} - \frac{16V^2d^3}{b^4\Gamma_4^2} + \frac{32cd^2V^2}{b^3\Gamma_4^2} - \frac{64edV^2}{b^2\Gamma_4^2} - 256\alpha_1^2 +$$
$$1024\alpha_2 + \frac{128V^2\left(g - \frac{f\alpha_1}{2} + \frac{d^2\alpha_1}{8b}\right)}{b\Gamma_4^2} = 0$$

(84)

We use equations (eq.75) and (eq.78) to calculate the value of $\alpha_4$ and then express $\alpha_2$ as shown in (eq.85) where we rely on replacing $\frac{\Gamma_4}{\alpha_3}$ with its constant value $V$, which is presented in (eq.80).

$$1024\alpha_2 = \frac{1024\Gamma^2}{V} - \frac{512d}{b} - \frac{6144\Gamma^4}{V^4b^4} + \frac{3072\Gamma^4}{V^3b^3} + \frac{384\Gamma^4}{V^2b^2} - \frac{4608d\Gamma^2}{V^2b^4} - \frac{1536\Gamma^2\alpha_1}{V^2b^2} + \frac{3072c\Gamma^2}{V^2b^3} + \frac{384f}{b^3} +$$
$$\frac{96V\left(f - \frac{d^2}{4b}\right)}{b^2} - \frac{672d^2}{b^4} + \frac{768cd}{b^3} - \frac{768e}{b^2} - \frac{192\Gamma^4}{Vb} + \frac{96fV}{b^2} + \frac{288d\Gamma^2}{b^2} + 96\Gamma^2\alpha_1 - \frac{192c\Gamma^2}{b} + \frac{24V^2\left(f - \frac{d^2}{4b}\right)}{b} -$$
$$\frac{168(d^2V)}{b^3} + \frac{192(cdV)}{b^2} - \frac{192eV}{b}$$

(85)

In order to pass from equation (eq.84) to equation (eq.86), we replace $\alpha_2$ with its shown expression in (eq.85).

$$4\Gamma_4^4 + \frac{1024\Gamma_4^4}{V} - \frac{512d}{b} - \frac{6144\Gamma_4^4}{V^4b^4} + \frac{4608\Gamma_4^4}{V^3b^3} - \frac{128\Gamma_4^4}{Vb} + \frac{448\Gamma_4^4}{V^2b^2} - \frac{5376d\Gamma_4^2}{V^2b^4} - \frac{2304\Gamma_4^2\alpha_1}{V^2b^2} + \frac{3584c\Gamma_4^2}{V^2b^3} + \frac{448f}{b^3} +$$
$$\frac{24V^2\left(f - \frac{d^2}{4b}\right)}{b} + \frac{192fV}{b^2} - \frac{192d^2V}{b^3} + \frac{192cdV}{b^2} - \frac{192eV}{b} + \frac{240d\Gamma_4^2}{b^2} + 48\alpha_1\Gamma_4^2 - \frac{160c\Gamma_4^2}{b} + \frac{384d\Gamma_4^2}{Vb^3} - \frac{128\Gamma_4^2\alpha_1}{Vb} -$$
$$\frac{256c\Gamma_4^2}{Vb^2} - \frac{880d^2}{b^4} + \frac{1024cd}{b^3} - \frac{1024e}{b^2} - \frac{16V^2d^3}{b^4\Gamma_4^2} + \frac{32cd^2V^2}{b^3\Gamma_4^2} - \frac{64edV^2}{b^2\Gamma_4^2} - 256\alpha_1^2 + \frac{128V^2g}{b\Gamma_4^2} - \frac{64V^2f\alpha_1}{b\Gamma_4^2} +$$
$$\frac{16V^2d^2\alpha_1}{b^2\Gamma_4^2} = 0$$

(86)

We replace $\alpha_1$ with its shown expression in (eq.82), and then we assemble terms of equation (eq.86) in function of degrees, in order to have the expression (eq.88) where coefficients are presented in (eq.54), (eq.55), (eq.56) and (eq.57).

$$\alpha_4 = x_1x_2x_3x_4 \Rightarrow 2048\alpha_4 = -32\Gamma_3^2\alpha_1 - \frac{8\Gamma_3^2\left(e - \frac{c^2}{4}\right)}{\alpha_3^2} + \frac{64\Gamma_3d}{\alpha_3} \quad (87)$$

$$\lambda_3(\Gamma_4^2)^4 + \lambda_2(\Gamma_4^2)^3 + \lambda_1(\Gamma_4^2)^2 + \lambda_0(\Gamma_4^2) = 0 \quad (88)$$

Since we supposed that $\left(\frac{4\Gamma_4^2\left(f-\frac{d^2}{4b}\right)}{b\alpha_3^2} + \frac{32f\Gamma_4}{b^2\alpha_3} + \frac{40(d^2\Gamma_4)}{b^3\alpha_3} - \frac{64cd\Gamma_4}{b^2\alpha_3} + \frac{64e\Gamma_4}{b\alpha_3}\right) = 0$ whereas adopting $\frac{\Gamma_4}{\alpha_3} \neq 0$, we eliminate the root zero as solution for polynomial equation (eq.88) and we use the cubic solution to solve the polynomial equation $(\lambda_3(\Gamma_4^2)^3 + \lambda_2(\Gamma_4^2)^2 + \lambda_1(\Gamma_4^2) + \lambda_0) = 0$, because all coefficients are expressed only in function of $c, d, e, f, g$ and the constant $V$. As a result, we have six possible values for $\Gamma_4$ as solutions for polynomial equation shown in (eq.53).

Supposing that $G_{\{\Gamma_4\}}$ is the group of solutions for shown equation in (eq.53), where these solutions are expressed as $\Gamma_{4,i}$ and $-\Gamma_{4,i}$ with $1 \leq i \leq 3$. The group of solutions $G_{\{\Gamma_4\}}$ is determined by relying on cubic root shown in (eq.90) and quadratic roots (eq.91) and (eq.92).

$$G_{\{\Gamma_4\}} = \{\Gamma_{4,1}, \Gamma_{4,2}, \Gamma_{4,3}, -\Gamma_{4,1}, -\Gamma_{4,2}, -\Gamma_{4,3}\} \quad (89)$$

We suppose that $b^i = \frac{\lambda_2}{\lambda_3}$, $c^i = \frac{\lambda_1}{\lambda_3}$ and $d^i = \frac{\lambda_0}{\lambda_3}$, whereas using the expressions (eq.54), (eq.55), (eq.56) and (eq.57). We suppose also that $D^i = 27d^i + 2b^{i3} - 9c^i b^i$ and $C^i = 9c^i - 3b^{i2}$. The solutions $\Gamma_{4,1}$, $\Gamma_{4,2}$ and $\Gamma_{4,3}$ for shown equation in (eq.53) are as follow:

$$\Gamma_{4,1}^2 = \frac{-b^i}{3} + \frac{1}{3}\sqrt[3]{-\frac{D^i}{2} + \sqrt{\left(\frac{D^i}{2}\right)^2 + \left(\frac{C^i}{3}\right)^3}} + \frac{1}{3}\sqrt[3]{-\frac{D^i}{2} - \sqrt{\left(\frac{D^i}{2}\right)^2 + \left(\frac{C^i}{3}\right)^3}} \quad (90)$$

$$\Gamma_{4,2}^2 = -\frac{\frac{b^i}{2} + \Gamma_{4,1}^2}{2} + \sqrt{\left(\frac{\frac{b^i}{2} + \Gamma_{4,1}^2}{2}\right)^2 - \frac{D^i}{\Gamma_{4,1}^2}} \quad (91)$$

$$\Gamma_{4,3}^2 = -\frac{\frac{b^i}{2} + \Gamma_{4,1}^2}{2} - \sqrt{\left(\frac{\frac{b^i}{2} + \Gamma_{4,1}^2}{2}\right)^2 - \frac{D^i}{\Gamma_{4,1}^2}} \quad (92)$$

After determining the values of $\Gamma_4$ by using cubic solution and quadratic solutions, the following step consists of solving the polynomial equation shown in (eq.72).

The coefficients $P$ in (eq.74) and $R$ in (eq.76) are dependent on $\Gamma_4^2$ and $\alpha_1$. The coefficient $\alpha_1$ in (eq.82) is dependent on $\Gamma_4^2$, whereas the coefficient $Q$ shown in (eq.75) is dependent on $\Gamma_4$ and $\alpha_1$. Therefore, we are going to use only $\Gamma_{4,1}$, $\Gamma_{4,2}$ and $\Gamma_{4,3}$ to calculate the potential values of $z$, because $\{-\Gamma_{4,1}, -\Gamma_{4,2}, -\Gamma_{4,3}\}$ are going only to inverse the sign of coefficient $Q$ and thereby inversing the signs of potential values of $z$ as solutions for polynomial equation shown in (eq.50), which will not influence the potential values of $x$ as solutions for sixth degree polynomial equation shown in (eq.48) because $x = \frac{1}{2}(z^2 - \alpha_1)$.

We use the first proposed theorem (Theorem 1) to calculate the four solutions for fourth degree polynomial equation shown in (eq.72) after calculating the coefficients $P$, $Q$ and $R$ for each value of $\Gamma_4$ from the group $\{\Gamma_{4,1}, \Gamma_{4,2}, \Gamma_{4,3}\}$. Thereby, we have twelve values to calculate as potential solutions for the polynomial equation shown in (eq.72).

After using first proposed theorem on (eq.72) for each value of $\Gamma_4$ from the group $\{\Gamma_{4,1}, \Gamma_{4,2}, \Gamma_{4,3}\}$, we have three groups of potential solutions for polynomial equation shown in (eq.72), where each group is dependent on different value of $\Gamma_4$. We express these groups of solutions as follow: $K_{\{\Gamma_{4,1}\}}, K_{\{\Gamma_{4,2}\}}, K_{\{\Gamma_{4,3}\}}$.

$$K_{\{\Gamma_{4,1}\}} = \left\{S_{(\Gamma_{4,1},1)}, S_{(\Gamma_{4,1},2)}, S_{(\Gamma_{4,1},3)}, S_{(\Gamma_{4,1},4)}\right\} \quad (93)$$

$$K_{\{\Gamma_{4,2}\}} = \left\{S_{(\Gamma_{4,2},1)}, S_{(\Gamma_{4,2},2)}, S_{(\Gamma_{4,2},3)}, S_{(\Gamma_{4,2},4)}\right\} \quad (94)$$

$$K_{\{\Gamma_{4,3}\}} = \left\{ S_{(\Gamma_{4,3},1)}, S_{(\Gamma_{4,3},2)}, S_{(\Gamma_{4,3},3)}, S_{(\Gamma_{4,3},4)} \right\} \quad (95)$$

Concerning the fourth degree polynomial equation shown in (eq.50), we have three groups of solutions where each group is dependent on different value of $\Gamma_4$; as shown in (eq.96), (eq.97) and (eq.98). The values of $S_{(\Gamma_{4,i},j)}$, where $1 \leq i \leq 3$ and $1 \leq j \leq 4$, are from the expressed solutions in the groups (eq.93), (eq.94) and (eq.95).

$$M_{\{\Gamma_{4,1}\}} = \left\{ -\frac{1}{4}\left[\frac{4\Gamma_{4,1}}{V} + \Gamma_{4,1}\right] + \frac{1}{4}S_{(\Gamma_{4,1},1)}, -\frac{1}{4}\left[\frac{4\Gamma_{4,1}}{V} + \Gamma_{4,1}\right] + \frac{1}{4}S_{(\Gamma_{4,1},2)}, -\frac{1}{4}\left[\frac{4\Gamma_{4,1}}{V} + \Gamma_{4,1}\right] + \right.$$
$$\left. \frac{1}{4}S_{(\Gamma_{4,1},3)}, -\frac{1}{4}\left[\frac{4\Gamma_{4,1}}{V} + \Gamma_{4,1}\right] + \frac{1}{4}S_{(\Gamma_{4,1},4)} \right\}$$

$$(96)$$

$$M_{\{\Gamma_{4,2}\}} = \left\{ -\frac{1}{4}\left[\frac{4\Gamma_{4,2}}{V} + \Gamma_{4,2}\right] + \frac{1}{4}S_{(\Gamma_{4,2},1)}, -\frac{1}{4}\left[\frac{4\Gamma_{4,2}}{V} + \Gamma_{4,2}\right] + \frac{1}{4}S_{(\Gamma_{4,2},2)}, -\frac{1}{4}\left[\frac{4\Gamma_{4,2}}{V} + \Gamma_{4,2}\right] + \right.$$
$$\left. \frac{1}{4}S_{(\Gamma_{4,2},3)}, -\frac{1}{4}\left[\frac{4\Gamma_{4,2}}{V} + \Gamma_{4,2}\right] + \frac{1}{4}S_{(\Gamma_{4,2},4)} \right\}$$

$$(97)$$

$$M_{\{\Gamma_{4,3}\}} = \left\{ -\frac{1}{4}\left[\frac{4\Gamma_{4,3}}{V} + \Gamma_{4,3}\right] + \frac{1}{4}S_{(\Gamma_{4,3},1)}, -\frac{1}{4}\left[\frac{4\Gamma_{4,3}}{V} + \Gamma_{4,3}\right] + \frac{1}{4}S_{(\Gamma_{4,3},2)}, -\frac{1}{4}\left[\frac{4\Gamma_{4,3}}{V} + \Gamma_{4,3}\right] + \right.$$
$$\left. \frac{1}{4}S_{(\Gamma_{4,3},3)}, -\frac{1}{4}\left[\frac{4\Gamma_{4,3}}{V} + \Gamma_{4,3}\right] + \frac{1}{4}S_{(\Gamma_{4,3},4)} \right\}$$

$$(98)$$

We suppose that $\dot{S}_{(\Gamma_{4,i},j)} = \left( -\frac{1}{4}\left[\frac{4\Gamma_{4,i}}{V} + \Gamma_{4,i}\right] + \frac{1}{4}S_{(\Gamma_{4,i},j)} \right)$ where $1 \leq i \leq 3$ and $1 \leq j \leq 4$, in order to simplify the expressed values in (eq.96), (eq.97) and (eq.98). Thereby, we have three groups of values as potential solutions for sixth degree polynomial equation shown in (eq.48). These three groups are as shown in (eq.99), (eq.100) and (eq.101) where $\alpha_{(1,\Gamma_{4,i})}$ is as follow:

$$\alpha_{(1,\Gamma_{4,i})} = \frac{\dfrac{\Gamma_{4,i}^4 + \dfrac{32\Gamma_{4,i}^4}{V^2 b^2} - \dfrac{8\Gamma_{4,i}^4}{Vb} + \dfrac{12d\Gamma_{4,i}^2}{b^2} - \dfrac{8c\Gamma_{4,i}^2}{b} - \dfrac{V^2\left(f - \dfrac{d^2}{4b}\right)}{b}}{4\Gamma_{4,i}^2}}{}$$

The expression of $\alpha_{(1,\Gamma_{4,i})}$ is an extending of the shown expression of $\alpha_1$ in (eq.82) by changing the value of $\Gamma_4$, where $\Gamma_{4,i}$ belong to the group $\left\{\Gamma_{4,1}, \Gamma_{4,2}, \Gamma_{4,3}\right\}$.

$$N_{\{\Gamma_{4,1}\}} = \left\{ \frac{1}{2}\left[\dot{S}_{(\Gamma_{4,1},1)}^2 - \alpha_{(1,\Gamma_{4,1})}\right], \frac{1}{2}\left[\dot{S}_{(\Gamma_{4,1},2)}^2 - \alpha_{(1,\Gamma_{4,1})}\right], \frac{1}{2}\left[\dot{S}_{(\Gamma_{4,1},3)}^2 - \alpha_{(1,\Gamma_{4,1})}\right], \frac{1}{2}\left[\dot{S}_{(\Gamma_{4,1},4)}^2 - \right.\right.$$
$$\left.\left. \alpha_{(1,\Gamma_{4,1})}\right]\right\}$$

$$(99)$$

$$N_{\{\Gamma_{4,2}\}} = \left\{ \frac{1}{2}\left[\dot{S}_{(\Gamma_{4,2},1)}^2 - \alpha_{(1,\Gamma_{4,2})}\right], \frac{1}{2}\left[\dot{S}_{(\Gamma_{4,2},2)}^2 - \alpha_{(1,\Gamma_{4,2})}\right], \frac{1}{2}\left[\dot{S}_{(\Gamma_{4,2},3)}^2 - \alpha_{(1,\Gamma_{4,2})}\right], \frac{1}{2}\left[\dot{S}_{(\Gamma_{4,2},4)}^2 - \right.\right.$$
$$\left.\left. \alpha_{(1,\Gamma_{4,2})}\right]\right\}$$

$$(100)$$

$$N_{\{\Gamma_{4,3}\}} = \left\{ \frac{1}{2}\left[\dot{S}_{(\Gamma_{4,3},1)}^2 - \alpha_{(1,\Gamma_{4,3})}\right], \frac{1}{2}\left[\dot{S}_{(\Gamma_{4,3},2)}^2 - \alpha_{(1,\Gamma_{4,3})}\right], \frac{1}{2}\left[\dot{S}_{(\Gamma_{4,3},3)}^2 - \alpha_{(1,\Gamma_{4,3})}\right], \frac{1}{2}\left[\dot{S}_{(\Gamma_{4,3},4)}^2 - \right.\right.$$
$$\left.\left. \alpha_{(1,\Gamma_{4,3})}\right]\right\}$$

$$(101)$$

We have twelve values as potential solutions for the sixth degree polynomial equation shown in (eq.48). However, many of them are only redundancies of others and there are only six official solutions to determine. The variables $\left\{\Gamma_{4,1}, \Gamma_{4,2}, \Gamma_{4,3}\right\}$ are the responsible of solution redundancies from one group to other.

In order to avoid the complications of calculating the twelve values from the groups $N_{\{\Gamma_{4,1}\}}$, $N_{\{\Gamma_{4,2}\}}$ and $N_{\{\Gamma_{4,3}\}}$, and then determining the six solutions for sixth degree polynomial equation shown in (eq.48), we propose the expressed values in (eq.102), (eq.103), (eq.104), (eq.105), (eq.106) and (eq.107) as the six official solutions for sixth degree polynomial equation shown in (eq.48).

The first four proposed values as solutions for sixth degree polynomial equation shown in (eq.48) are from the group $N_{\{\Gamma_{4,1}\}}$, whereas the fifth and sixth values are expressed by deduction using the expressions of quadratic roots. The useless redundancies of solutions are from one group to other; therefore, we choose the first four solutions from the same group $N_{\{\Gamma_{4,1}\}}$.

In our six proposed solutions, we use the values $\dot{S}_{(\Gamma_{4,1},j)} = -\frac{1}{4}\Gamma_{4,1} + \frac{1}{4}S_{(\Gamma_{4,1},j)}$ where $1 \leq j \leq 4$ and $S_{(\Gamma_{4,1},j)}$ from $K_{\{\Gamma_{4,1}\}}$. The variable $\alpha_{(1,\Gamma_{4,i})}$ is expressed as follow:

$$\alpha_{(1,\Gamma_{4,i})} = \frac{\Gamma_{4,i}^4 + \frac{32\Gamma_{4,i}^4}{v^2 b^2} - \frac{8\Gamma_{4,i}^4}{vb} + \frac{12d\Gamma_{4,i}^2}{b^2} - \frac{8c\Gamma_{4,i}^2}{b} - \frac{v^2\left(f - \frac{d^2}{4b}\right)}{b}}{4\Gamma_{4,i}^2}$$

$\Gamma_{4,1}$ is from the group $G_{\{\Gamma_4\}}$ shown in (eq.88) which contains the solutions for polynomial equation (eq.52). The values of $\dot{S}_{(\Gamma_{4,1},1)}$, where $1 \leq j \leq 4$, are the solutions for quartic equation (eq.50) and they are determined by using Theorem 1.

$$S_1 = \frac{1}{2}\left[\dot{S}_{(\Gamma_{4,1},1)}^2 - \alpha_{(1,\Gamma_{4,1})}\right] \quad (102)$$

$$S_2 = \frac{1}{2}\left[\dot{S}_{(\Gamma_{4,1},2)}^2 - \alpha_{(1,\Gamma_{4,1})}\right] \quad (103)$$

$$S_3 = \frac{1}{2}\left[\dot{S}_{(\Gamma_{4,1},3)}^2 - \alpha_{(1,\Gamma_{4,1})}\right] \quad (104)$$

$$S_4 = \frac{1}{2}\left[\dot{S}_{(\Gamma_{4,1},4)}^2 - \alpha_{(1,\Gamma_{4,1})}\right] \quad (105)$$

$$S_5 = -\frac{b - S_1 - S_2 - S_3 - S_4}{2} - \sqrt{\left(\frac{b - S_1 - S_2 - S_3 - S_4}{2}\right)^2 - \frac{g}{S_1 S_2 S_3 S_4}} \quad (106)$$

$$S_6 = -\frac{b - S_1 - S_2 - S_3 - S_4}{2} + \sqrt{\left(\frac{b - S_1 - S_2 - S_3 - S_4}{2}\right)^2 - \frac{g}{S_1 S_2 S_3 S_4}} \quad (107)$$

### 3.3. Third Proposed Theorem

In the second proposed theorem in this paper, we rely on reducing sixth degree polynomial equation to fourth degree in order to find proper roots for concerned equation. However, the fourth degree part of resulted quartic equation by reduction is dependent on the fifth degree part of concerned sixth degree polynomial. Thereby, the equation $Aw^6 + Cw^4 + Dw^3 + Ew^2 + Fw + G = 0$ with $A \neq 0$ where the coefficient of fifth degree part is equal zero imposes a problem of reduction to fourth degree.

Therefore, in this subsection, we present a third theorem to propose six formulary solutions for sixth degree polynomial equation in general form $Aw^6 + Cw^4 + Dw^3 + Ew^2 + Fw + G = 0$ with $A \neq 0$, where coefficients belong to the group of numbers $\mathbb{R}$ and the coefficient of fifth degree part equal zero.

This third proposed theorem is based on the same logic and calculations of Theorem 2. However, it is distinguished by treating sixth degree polynomial where the fifth degree part is absent. First, we pass from equation expression $Aw^6 + Cw^4 + Dw^3 + Ew^2 + Fw + G = 0$ with $A \neq 0$ to the expression $w^6 + \frac{Cw^4}{A} + \frac{Dw^3}{A} + \frac{Ew^2}{A} + \frac{Fw}{A} + \frac{G}{A} = 0$, and then we use the expression $w = \sqrt{\frac{-C}{15A}} + x$ to induce a fifth degree part whereas eliminating the fourth degree part of concerned sixth degree polynomial.

Inducing the fifth degree part in polynomial equation shown in (eq.108), whereas eliminating the fourth degree part; help to reduce the amount of calculations comparing to the situation of having a fourth degree element to treat during the process of reduction from sixth degree polynomial equation to quartic equation.

The shown equation in (eq.108) is resulted by replacing $w$ with $w = \sqrt{\frac{-C}{15A}} + x$. The coefficients of equation (eq.108) are as expressed in (eq.109), (eq.110), (eq.111), (eq.112) and (eq.113).

$$x^6 + bx^5 + dx^3 + ex^2 + fx + g = 0 \quad (108)$$

$$b = 6\sqrt{\frac{-C}{15A}} \quad (109)$$

$$d = \frac{8C}{3A}\sqrt{\frac{-C}{15A}} + \frac{D}{A} \quad (110)$$

$$e = \frac{-C^2}{3A^2} + \frac{3D}{A}\sqrt{\frac{-C}{15A}} + \frac{E}{A} \quad (111)$$

$$f = -\frac{18C^2}{5A^2}\sqrt{\frac{-C}{15A}} - \frac{DC}{5A^2} + \frac{2E}{A}\sqrt{\frac{-C}{15A}} + \frac{F}{A} \quad (112)$$

$$g = \frac{-16C^3}{3375A^3} - \frac{DC}{15A^2}\sqrt{\frac{-C}{15A}} - \frac{EC}{15A^2} + \frac{F}{A}\sqrt{\frac{-C}{15A}} + \frac{G}{A} \quad (113)$$

***Theorem 3***

In order to reduce the sixth degree polynomial equation $Aw^6 + Cw^4 + Dw^3 + Ew^2 + Fw + G = 0 \text{ with } A \neq 0$ to the quartic equation shown in (eq.114), where coefficients belong to the group of numbers $\mathbb{R}$ we first replace $w$ with $w = \sqrt{\frac{-C}{15A}} + x$. Then, the reduction from sixth degree to fourth degree is conducted by supposing $x = (x_0 x_1 + x_0 x_2 + x_0 x_3 + x_1 x_2 + x_1 x_3 + x_2 x_3)$, whereas supposing $z = (x_0 + x_1 + x_2 + x_3)$ is the solution for fourth degree polynomial equation in (eq.114) by using Theorem 1 and relying on the expression $x_3 = -\frac{Y_3}{4}$. The variable $Y_3$ is defined as shown in (eq.115) where $\alpha_3$ is presented in (eq.119) and $Y_4$ is the solution for the polynomial equation (eq.120), which relies on the coefficients (eq.121), (eq.122), (eq.123) and (eq.124). The shown coefficients in (eq.121), (eq.122), (eq.123) and (eq.124) are expressed by using the constant $V$, which is defined in (eq.125). The coefficients $Y_3, Y_2, Y_1$ and $Y_0$ of quartic equation (eq.114) are determined by using calculated value of $Y_4$ and using the shown expressions in (eq.115), (eq.116), (eq.117) and (eq.118). The six proposed solutions for polynomial equation $Aw^6 + Cw^4 + Dw^3 + Ew^2 + Fw + G = 0 \text{ with } A \neq 0$ are as shown in (eq.136), (eq.137), (eq.138), (eq.139), (140) and (eq.141).

$$z^4 + Y_3 z^3 + Y_2 z^2 + Y_1 z + Y_0 = 0 \quad (114)$$

$$Y_3 = \frac{4\alpha_3}{b} + Y_4 \quad (115)$$

$$Y_2 = \frac{8Y_4^2}{Vb} - \frac{6d}{b^2} + \frac{\left(f - \frac{d^2}{4b}\right)V^2}{2bY_4^2} - \frac{8Y_4^2}{V^2 b^2} \quad (116)$$

$$\Upsilon_1 = \frac{5\Upsilon_4^3}{vb} + \frac{3Vd^2}{4b^3\Upsilon_4} - \frac{6d\Upsilon_4}{b^2} + \frac{eV}{b\Upsilon_4} - \frac{\Upsilon_4^3}{4} - \frac{8\Upsilon_4^3}{v^2b^2} + \frac{f - \frac{d^2}{4b}}{4\Upsilon_4 b}V^2 \quad (117)$$

$$\Upsilon_0 = \frac{\Upsilon_4^4}{2Vb} - \frac{V^2d^3}{16b^4\Upsilon_4^2} + \frac{3Vd^2}{8b^3} - \frac{3d\Upsilon_4^2}{4b^2} + \frac{eV}{2b} - \frac{eV^2d}{4b^2\Upsilon_4^2} + \frac{gV^2}{2b\Upsilon_4^2} - \left(\frac{\Upsilon_4^2}{4} + V^2\frac{f - \frac{d^2}{4b}}{4b\Upsilon_4^2}\right)\left(\frac{\Upsilon_4^2}{4} + \frac{8\Upsilon_4^3}{V^2b^2} - \frac{2\Upsilon_4^2}{VB} + \frac{3d}{b^2} - \frac{\left(f - \frac{d^2}{4b}\right)}{4b\Upsilon_4^2}V^2\right)$$

$$(118)$$

$$\alpha_3 = -\frac{\Upsilon_4\frac{4\left(f - \frac{d^2}{4b}\right)}{b}}{\frac{32f}{b^2} + \frac{40d^2}{b^3} + \frac{64e}{b}} \quad (119)$$

$$\beta_3\Upsilon_4^6 + \beta_2\Upsilon_4^4 + \beta_1\Upsilon_4^2 + \beta_0 = 0 \quad (120)$$

$$\beta_3 = -\frac{40960}{v^4b^4} + \frac{16384}{v^3b^3} - \frac{1536}{v^2b^2} \quad (121)$$

$$\beta_2 = -\frac{24576d}{v^2b^4} + \frac{3072d}{vb^3} + \frac{1024}{v} \quad (122)$$

$$\beta_1 = -\frac{512d}{b} + \frac{1536f}{b^3} + \frac{28V^2f}{b} - \frac{7V^2d^2}{b^2} + \frac{96Vf}{b^2} - \frac{168d^2V}{b^3} - \frac{192Ve}{b} - \frac{3456d^2}{b^4} - \frac{1024e}{b^2}$$

$$(123)$$

$$\beta_0 = -\frac{64V^2d^3}{b^4} - \frac{64eV^2d}{b^2} + \frac{128V^2g}{b} + \frac{192V^2df}{b^3} \quad (124)$$

$$V = -\frac{\frac{32f}{b^2} + \frac{40d^2}{b^3} + \frac{64e}{b}}{\frac{4\left(f - \frac{d^2}{4b}\right)}{b}} \quad (125)$$

### 3.4. *Proof of Theorem 3*

Considering the sixth degree polynomial equation $Aw^6 + Cw^4 + Dw^3 + Ew^2 + Fw + G = 0$ where the fifth degree part is absent, we divide the polynomial on $A$ and then we use the expression $w = \sqrt{\frac{-C}{15A}} + x$ in order to induce a fifth degree part and eliminate the fourth degree part. Then, by using the expression (eq.62), we reduce the resulted sixth degree polynomial (eq.108) to the quartic polynomial shown in (eq.126).

$$v_4z^4 + v_3z^3 + v_2z^2 + v_1z + v_0 = 0 \quad (126)$$

We rely on the expressions of $\{\alpha_1, \alpha_2, \alpha_3, \alpha_4\}$ in (eq.64), $x$ in (eq.65), $x^2$ in (eq.66), $x^3$ in (eq.67), $x^5$ in (eq.69) and $x^6$ in (eq.70) to express the fourth degree polynomial shown in (eq.126) where the values of coefficients are as follow:

$$v_4 = 2b\alpha_3^2$$

$$v_3 = 8\alpha_3^3 + d\alpha_3 + 2b[\alpha_2 + 6\alpha_4]\alpha_3$$

$$v_2 = 12[\alpha_2 + 6\alpha_4]\alpha_3^2 + \frac{1}{2}b[(\alpha_2 + 6\alpha_4)^2 - 4\alpha_3^2\alpha_1] + \frac{1}{2}d[\alpha_2 + 6\alpha_4] + \frac{1}{2}f$$

$$v_1 = 6[(\alpha_2 + 6\alpha_4)]^2\alpha_3 - 2b\alpha_3\alpha_1[\alpha_2 + 6\alpha_4] - d\alpha_1\alpha_3 + 2e\alpha_3$$

$$v_0 = [\alpha_2 + 6\alpha_4]^3 - \frac{1}{2}b\alpha_1[\alpha_2 + 6\alpha_4]^2 - \frac{1}{2}d[\alpha_2 + 6\alpha_4]\alpha_1 + e[\alpha_2 + 6\alpha_4] - \frac{1}{2}f\alpha_1 + g$$

We have the shown fourth degree polynomial in (eq.114) after dividing the polynomial (eq.126) on $v_4$. The coefficients of polynomial (eq.114) are as follow:

$$Y_4 = \frac{[\alpha_2 + 6\alpha_4] + \frac{d}{2b}}{\alpha_3} \Rightarrow \alpha_2 = \alpha_3 Y_4 - \frac{d}{2b} - 6\alpha_4$$

$$Y_3 = \frac{4\alpha_3}{b} + Y_4$$

$$Y_2 = \frac{6\left[Y_4\alpha_3 - \frac{d}{2b}\right]}{b} + \frac{Y_4^2}{4} - \alpha_1 + \frac{\left(f - \frac{d^2}{4b}\right)}{4b\alpha_3^2}$$

$$Y_1 = \frac{3\left[Y_4\alpha_3 - \frac{d}{2b}\right]^2}{b\alpha_3} - Y_4\alpha_1 + \frac{e}{b\alpha_3}$$

$$Y_0 = \frac{\left[Y_4\alpha_3 - \frac{d}{2b}\right]^3}{2b\alpha_3^2} - \frac{1}{4}\alpha_1 Y_4^2 + \frac{e\left[\alpha_3 Y_4 - \frac{d}{2b}\right]}{2b\alpha_3^2} + \frac{g - \frac{f\alpha_1}{2} + \frac{d^2\alpha_1}{8b}}{2b\alpha_3^2}$$

We have the fourth degree polynomial $y^4 + My^2 + Ny + O = 0$ by replacing $z$ with $\frac{-Y_3 + y}{4}$ in the polynomial (eq.114). The coefficients $M$, $N$ and $O$ are as expressed in (eq.127), (eq.128) and (eq.129).

$$M = -2Y_4^2 - \frac{96\alpha_3^2}{b^2} + \frac{48Y_4\alpha_3}{b} - \frac{48d}{b^2} - 16\alpha_1 + \frac{4\left(f - \frac{d^2}{4b}\right)}{b\alpha_3^2} \quad (127)$$

$$N = \frac{512\alpha_3^3}{b^3} - \frac{384\alpha_3^2 Y_4}{b^2} + \frac{64\alpha_3 Y_4^2}{b} + \frac{384d\alpha_3}{b^3} + \frac{128\alpha_3\alpha_1}{b} - \frac{32f}{b^2\alpha_3} - \frac{96dY_4}{b^2} - 32Y_4\alpha_1 - \frac{8Y_4\left(f - \frac{d^2}{4b}\right)}{b\alpha_3^2} + \frac{56d^2}{b^3\alpha_3} + \frac{64e}{b\alpha_3}$$

(128)

$$O = -\frac{768\alpha_3^4}{b^4} + Y_4^4 + \frac{768Y_4\alpha_3^3}{b^3} + \frac{16Y_4^3\alpha_3}{b} - \frac{224Y_4^2\alpha_3^2}{b^2} - \frac{768d\alpha_3^2}{b^4} - \frac{256\alpha_3^2\alpha_1}{b^2} + \frac{64f}{b^3} - \frac{48dY_4^2}{b^2} - 16\alpha_1 Y_4^2 + \frac{4\left(f - \frac{d^2}{4b}\right)}{b\alpha_3^2}Y_4^2 + \frac{384dY_4\alpha_3}{b^3} + \frac{128Y_4\alpha_3\alpha_1}{b} + \frac{32fY_4}{b^2\alpha_3} + \frac{40d^2 Y_4}{b^3\alpha_3} - \frac{208d^2}{b^4} - \frac{256e}{b^2} + \frac{64eY_4}{b\alpha_3} - \frac{16d^3}{b^4\alpha_3^2} - \frac{64ed}{b^2\alpha_3^2} + \frac{128\left(g - \frac{f\alpha_1}{2} + \frac{d^2\alpha_1}{8b}\right)}{b\alpha_3^2}$$

(129)

In order to reduce the expression of $O$ in (eq.129), and find a way to determine the value of $Y_4$, we suppose that $\left(\frac{4Y_4^2\left(f - \frac{d^2}{4b}\right)}{b\alpha_3^2} + \frac{32fY_4}{b^2\alpha_3} + \frac{40d^2 Y_4}{b^3\alpha_3} + \frac{64eY_4}{b\alpha_3}\right) = 0$ where $\frac{Y_4}{\alpha_3} \neq 0$. As a result, we have the shown expression in (eq.130).

$$\frac{Y_4}{\alpha_3} = V = -\frac{\left(\frac{32f}{b^2} + \frac{40d^2}{b^3} + \frac{64e}{b}\right)}{\frac{4\left(f - \frac{d^2}{4b}\right)}{b}} \quad (130)$$

We rely on the same used logic and processes of calculation in the proof of Theorem 2, in order to continue the proof of third proposed theorem. Thereby, the variables $\alpha_1$, $\alpha_4$ and $\alpha_2$ are as expressed in (eq.131), (eq.132) and (eq.133) respectively, whereas the resulted polynomial equation to determine the value of $Y_4$ is as shown in (eq.134).

$$\alpha_1 = \frac{Y_4^4 + \frac{32Y_4^4}{V^2 b^2} - \frac{8Y_4^4}{Vb} + \frac{12dY_4^2}{b^2} - \frac{V^2\left(f - \frac{d^2}{4b}\right)}{b}}{4Y_4^2} \quad (131)$$

$$2048\alpha_4 = \frac{2048Y^4}{V^4 b^4} - \frac{1024Y^4}{V^3 b^3} - \frac{128Y^4}{V^2 b^2} + \frac{1536dY^2}{V^2 b^4} + \frac{512Y^2\alpha_1}{V^2 b^2} - \frac{128f}{b^3} - \frac{32V\left(f - \frac{d^2}{4b}\right)}{b^2} + \frac{224d^2}{b^4} + \frac{256e}{b^2} + \frac{64Y^4}{Vb} -$$
$$\frac{32fV}{b^2} - \frac{96dY^2}{b^2} - 32Y^2\alpha_1 - \frac{8V^2\left(f - \frac{d^2}{4b}\right)}{b} + \frac{56d^2V}{b^3} + \frac{64eV}{b}$$

(132)

$$1024\alpha_2 = \frac{1024Y^2}{V} - \frac{512d}{b} - \frac{6144Y^4}{V^4 b^4} + \frac{3072Y^4}{V^3 b^3} + \frac{384Y^4}{V^2 b^2} - \frac{4608dY^2}{V^2 b^4} - \frac{1536Y^2\alpha_1}{V^2 b^2} + \frac{384f}{b^3} + \frac{96V\left(f - \frac{d^2}{4b}\right)}{b^2} -$$
$$\frac{672d^2}{b^4} - \frac{768e}{b^2} - \frac{192Y^4}{Vb} + \frac{96fV}{b^2} + \frac{288dY^2}{b^2} + 96Y^2\alpha_1 + \frac{24V^2\left(f - \frac{d^2}{4b}\right)}{b} - \frac{168d^2V}{b^3} - \frac{192eV}{b}$$

(133)

$$4Y_4^4 + \frac{1024Y_4^2}{V} - \frac{512d}{b} - \frac{6144Y_4^4}{V^4 b^4} + \frac{4608Y_4^4}{V^3 b^3} - \frac{128Y_4^4}{V^2 b^2} + \frac{448Y_4^4}{Vb} - \frac{5376dY_4^2}{V^2 b^2} - \frac{2304Y_4^2\alpha_1}{V^2 b^2} + \frac{448f}{b^3} + \frac{24V^2\left(f - \frac{d^2}{4b}\right)}{b} +$$
$$\frac{192fV}{b^2} - \frac{192d^2V}{b^3} - \frac{192eV}{b} + \frac{240dY_4^2}{b^2} + 48\alpha_1 Y_4^2 + \frac{384dY_4^2}{Vb^3} - \frac{128Y_4^2\alpha_1}{Vb} - \frac{880d^2}{b^4} - \frac{1024e}{b^2} - \frac{16V^2 d^3}{b^4 Y_4^2} -$$
$$\frac{64edV^2}{b^2 Y_4^2} - 256\alpha_1^2 + \frac{128V^2 g}{bY_4^2} - \frac{64V^2 f\alpha_1}{bY_4^2} + \frac{16V^2 d^2\alpha_1}{b^2 Y_4^2} = 0$$

(134)

We use the shown value of $\frac{Y_4}{\alpha_3}$ in (eq.130) and we replace $\alpha_1$ and $\alpha_2$ with their shown expressions in (eq.131) and (eq.133), in order to pass from equation (eq.134) to polynomial expression $\beta_3(Y_4^2)^4 + \beta_2(Y_4^2)^3 + \beta_1(Y_4^2)^2 + \beta_0(Y_4^2) = 0$ where coefficients are as presented in (eq.121), (eq.122), (eq.123) and (eq.124).

Relying on the proof of Theorem 2, we calculate only the expression $Y_{4,1}$ shown in (eq.135) as a root for expressed equation in (eq.120), and then we determine the roots of quartic equation $y^4 + My^2 + Ny + O = 0$. Therefore, we start by calculating the values of $\alpha_1$, $\alpha_4$ and $\alpha_2$ by replacing the variable $Y_4$ with the value of $Y_{4,1}$, then we calculate the values of $M, N$ and $O$, and finally we finish by using Theorem 1.

We suppose that $b^i = \frac{\beta_2}{\beta_3}, c^i = \frac{\beta_1}{\beta_3}$ and $d^i = \frac{\beta_0}{\beta_3}$, whereas using the expressions (eq.121), (eq.122), (eq.123) and (eq.124). We suppose also that $D^i = 27d^i + 2b^{i^3} - 9c^i b^i$ and $C^i = 9c^i - 3b^{i^2}$. The solution $Y_{4,1}$ for shown equation in (eq.120) is as follow:

$$Y_{4,1}^2 = \frac{-b^i}{3} + \frac{1}{3}\sqrt[3]{-\frac{D^i}{2} + \sqrt{\left(\frac{D^i}{2}\right)^2 + \left(\frac{C^i}{3}\right)^3}} + \frac{1}{3}\sqrt[3]{-\frac{D^i}{2} - \sqrt{\left(\frac{D^i}{2}\right)^2 + \left(\frac{C^i}{3}\right)^3}} \quad (135)$$

As we mentioned in the proof of Theorem 2, the use of other roots of polynomial (eq.120) in the quartic polynomial $y^4 + My^2 + Ny + O = 0$ to calculate the values of $M, N$ and $O$ generates redundancies of roots for the sixth degree polynomial equation $Aw^6 + Cw^4 + Dw^3 + Ew^2 + Fw + G = 0$.

We determine the group $K^i_{\{Y_{4,1}\}}$, which contains the four roots for quartic equation $y^4 + My^2 + Ny + O = 0$, by using Theorem 1.

$$K^i_{\{Y_{4,1}\}} = \left\{ S_{(Y_{4,1},1)}, S_{(Y_{4,1},2)}, S_{(Y_{4,1},3)}, S_{(Y_{4,1},4)} \right\}$$

We present the group of roots for quartic equation shown in (eq.114) as $M^i_{\{Y_{4,1}\}}$, which is determined by relying on the group $K^i_{\{Y_{4,1}\}}$.

$$M^i_{\{Y_{4,1}\}} = \left\{ -\frac{1}{4}\left[\frac{4Y_{4,1}}{V} + Y_{4,1}\right] + \frac{1}{4}S_{(Y_{4,1},1)}, -\frac{1}{4}\left[\frac{4Y_{4,1}}{V} + Y_{4,1}\right] + \frac{1}{4}S_{(Y_{4,1},2)}, -\frac{1}{4}\left[\frac{4Y_{4,1}}{V} + Y_{4,1}\right] \right.$$
$$\left. + \frac{1}{4}S_{(Y_{4,1},3)}, -\frac{1}{4}\left[\frac{4Y_{4,1}}{V} + Y_{4,1}\right] + \frac{1}{4}S_{(Y_{4,1},4)} \right\}$$

We present each root for quartic equation shown in (eq.114) as $\xi_{(Y_{4,1},i)} = -\frac{1}{4}\left[\frac{4Y_{4,1}}{V} + Y_{4,1}\right] + \frac{1}{4}S_{(Y_{4,1},i)}$, whereas $S_{(Y_{4,1},i)}$ is from the group $K^i_{\{Y_{4,1}\}}$.

The proposed six solutions for sixth degree polynomial equation $Aw^6 + Cw^4 + Dw^3 + Ew^2 + Fw + G = 0$ are as expressed in (eq.136), (eq.137), (eq.138), (eq.139), (140) and (eq.141). The expressions $\xi_{(Y_{4,1},1)}, \xi_{(Y_{4,1},2)}, \xi_{(Y_{4,1},3)}$ and $\xi_{(Y_{4,1},4)}$ present the calculated roots for quartic equation (eq.114) by using Theorem 1. $Y_{4,1}$ is calculated by using the shown expression in (eq.135). We use the expression (eq.131) to calculate the value of $\alpha_{(1,Y_{3,1})}$; thereby, its value is as follow:

$$\alpha_{(1,Y_{3,1})} = \frac{Y_{4,1}^4 + \frac{32Y_{4,1}^4}{V^2 b^2} - \frac{8Y_{4,1}^4}{Vb} + \frac{12dY_{4,1}^2}{b^2} - \frac{V^2\left(f - \frac{d^2}{4b}\right)}{b}}{4Y_{4,1}^2}$$

$$S_1 = \frac{1}{2}\left[\xi_{(Y_{4,1},1)}^2 - \alpha_{(1,Y_{4,1})}\right] \quad (136)$$

$$S_2 = \frac{1}{2}\left[\xi_{(Y_{4,1},2)}^2 - \alpha_{(1,Y_{4,1})}\right] \quad (137)$$

$$S_3 = \frac{1}{2}\left[\xi_{(Y_{4,1},3)}^2 - \alpha_{(1,Y_{4,1})}\right] \quad (138)$$

$$S_4 = \frac{1}{2}\left[\xi_{(Y_{4,1},4)}^2 - \alpha_{(1,Y_{4,1})}\right] \quad (139)$$

$$S_5 = -\frac{b - S_1 - S_2 - S_3 - S_4}{2} - \sqrt{\left(\frac{b - S_1 - S_2 - S_3 - S_4}{2}\right)^2 - \frac{g}{S_1 S_2 S_3 S_4}} \quad (140)$$

$$S_6 = -\frac{b - S_1 - S_2 - S_3 - S_4}{2} + \sqrt{\left(\frac{b - S_1 - S_2 - S_3 - S_4}{2}\right)^2 - \frac{g}{S_1 S_2 S_3 S_4}} \quad (141)$$

## 4. Conclusion

In the first presented theorem, we propose new four formulary solutions for any quartic polynomial equation in general form, which enabled us to develop the formulary structures of six roots for any sixth degree polynomial equation in general form.

The proposed expressions as solutions in first theorem enable to calculate the four roots of any quartic polynomial equation nearly simultaneously, whereas the proposed solutions in second theorem and third theorem enable to calculate the six roots of sixth degree polynomial equations in general forms nearly simultaneously.

The third proposed theorem is based on the same logic and calculations of Theorem 2 whereas proposing six roots for sixth degree polynomials. However, it is distinguished by treating the specific form $Aw^6 + Cw^4 + Dw^3 + Ew^2 + Fw + G = 0$ $with\ A \neq 0$ where the fifth degree part is absent, which is essential to reduce sixth degree polynomial equation to quartic equation. The third theorem is also distinguished by eliminating the fourth degree part of concerned sixth degree polynomial, in order to reduce the amount of calculations.

The principal criteria of presented theorems is proposing new radical solutions for quartic equations and sixth degree polynomial equations to enable the calculation of all roots of these equations nearly simultaneously, whereas using the expressions of quadratic roots and cubic roots as subparts of each proposed solution.

## References

[1]  G.Cardano, Artis Magnae, Sive de Regulis Algebraicis Liber Unus, 1545. English transl.: The Great Art, or The Rules of Algebra.

Translated and edited by Witmer, T. R. MIT Press, Cambridge, Mass, 1968.

[2]    D.S.Dumit, R.M. Foote, Abstract algebra. John Wiley, pp. 606-616, 2004.

[3]    T.J.Osler, "Cardan polynomials and the reduction of radicals," *Mathematics Magazine*, vol. 47, no. 1, pp. 26-32, 2001.

[4]    S.Janson, Roots of polynomials of degrees 3 and 4, 2010.

[5]    J. L.Lagrange, Réflexions sur la résolution algébrique des équations, in: ouvres de Lagrange, J.A. Serret ed & Gauthier-Villars, Vol. 3, pp. 205-421, 1869.

[6]    W. M.Faucette, "A geometric interpretation of the solution of the general quartic polynomial," Amer. Math. Monthly, vol. 103, pp. 51-57, 1996.

[7]    D.René, The Geometry of Rene Descartes with a facsimile of the 1st edition, Courier Corporation, 2012.

[8]    L.Euler, De formis radicum aequationum cuiusque ordinis coniectatio, 1738. English transl.: A conjecture on the forms of the roots of equations. Translated by Bell, J. Cornell University, NY, USA, 2008.

[9]    J.Bewersdor, Algebra fur Einsteiger, Friedr. Vieweg Sohn Verlag. English transl.: Galois Theory for Beginners. A Historical Perspective, 2004. Translated by Kramer, D. American Mathematical Society (AMS), Providence, R.I., USA, 2006.

[10]  H.Helfgott, M.Helfgott, A modern vision of the work of cardano and ferrari on quartics. Convergence (MAA).

[11]  M.Rosen, "Niels hendrik abel and equations of the fiffth degree," American Mathematical Monthly, vol. 102, no. 6, pp. 495-505, 1995.

[12]  Garling, D. J. H. Galois Theory. Cambridge Univ. Press, Mass., USA, 1986.

[13]  Grillet, P. Abstract Algebra. 2nd ed. Springer, New York, USA, 2007.

[14]  B. L.van der Waerden, Algebra, Springer-Verlag, Berlin, Vol. 1, 3rd ed. 1966. English transl.: Algebra. Translated by Schulenberg J.R. and Blum, F. Springer-Verlag, New York, 1991.

[15]  S. L.Shmakov, "A universal method of solving quartic equations," *International Journal of Pure and Applied Mathematics*, vol. 71, no. 2, pp. 251-259, 2011.

[16]  A.Fathi, and Sharifan, N. "A classic new method to solve quartic equations," *Applied and Computational Mathematics,* vol. 2, no. 2, pp. 24-27, 2013.

[17]  F.T.Tehrani, "Solution to polynomial equations, a new approach," *Applied Mathematics*, vol. 11, no. 2, pp. 53-66, 2020.

[18]  Y. J.Nahon, "Method for solving polynomial equations," *Journal of Applied & Computational Mathematics,* vol. 7, no. 3, pp. 2-12, 2018.

[19]  E.A.Tschirnhaus, "method for removing all intermediate terms from a given equation," *THIS BULLETIN*, vol. 37, no. 1, pp 1-3. Original: Methodus auferendi omnes terminos intermedios ex data equatione, Acta Eruditorium vol. 2, pp. 204-207, 1683.

[20]  V.S.Adamchik, D.J.Jeffrey, "Polynomial transformations of tschirnhaus, bring and jerrard," *ACM SIGSAM Bulletin*, vol. 37, no. 3, pp. 90-94, 2003.

[21]  D.Lazard, Solving quintics in radicals. In: Olav Arnfinn Laudal, Ragni Piene, The Legacy of Niels Henrik Abel. 1869, pp. 207-225. Springer, Berlin, Heidelberg.

[22]  D.S.Dummit, "Solving solvable quintics," Mathematics of Computation, vol. 57, no. 1, pp. 387-401, 1991.

[23]  BB. Samuel, "On Solvability of Higher Degree Polynomial Equations," Journal of Applied Science and Innovations, vol. 1, no. 2, 2017.

[24]  BB. Samuel, "The General Quintic Equation, its Solution by Factorization into Cubic and Quadratic Factors," Journal of Applied Science and Innovations, vol. 1, no. 2, 2017.

[25]  P. Pesic, Abel's Proof: An Essay on the Sources and Meaning of Mathematical Unsolvability, MIT Press, 2004.

[26]  P. L. Wantzel, Démonstration de l'impossibilité de résoudre toutes les équations algébriques avec des radicaux, J. Math. Pures Appl., 1re série, vol. 4, pp. 57-65, 1845.